\documentclass[]{spie}  %>>> use for US letter paper
\usepackage{color}
\usepackage{textcomp}
\usepackage[mathcal]{euscript}
\usepackage{mathrsfs}
\usepackage{subcaption}
\usepackage{amsfonts,amssymb,amsbsy}
\usepackage{graphicx,empheq}
\usepackage{nomencl}
\usepackage{arydshln}
\usepackage{multirow}
\usepackage[tableposition=top]{caption}
\usepackage{multirow}
\usepackage{mathtools}
\usepackage{lmodern}
\usepackage{arydshln}
\def\block(#1,#2)#3{\multicolumn{#2}{c}{\multirow{#1}{*}{$ #3 $}}}
\usepackage[font=footnotesize]{caption}
    \usepackage{booktabs}
    \usepackage[dvipsnames]{xcolor}
    \usepackage{tikz}
\usepackage[customcolors,shade]{hf-tikz}
\newcommand*{\boxcolor}{orange}
\makeatletter
\renewcommand{\boxed}[1]{\textcolor{\boxcolor}{%
\tikz[baseline={([yshift=-1ex]current bounding box.center)}] \node [rectangle, minimum width=1ex,rounded corners,draw] {\normalcolor\m@th$\displaystyle#1$};}}
 \makeatother

\def\e1{{\varepsilon_{1}}}
\def\b1{{\beta_{1}}}
\def\bp3{{\beta_{3}}}
\def\ep3{{\varepsilon_{3}}}

\newcommand{\mb}{\mathbf}

\newcommand{\mc}{\mathcal}

\newtheorem{rmk}{Remark}[section]

\MHInternalSyntaxOn
\providecommand*\phantomword[3][c]{%
\mathchoice
{\MT_phantom_word:NNnn #1\displaystyle {#2}{#3}}%
{\MT_phantom_word:NNnn #1\textstyle {#2}{#3}}%
{\MT_phantom_word:NNnn #1\scriptstyle {#2}{#3}}%
{\MT_phantom_word:NNnn #1\scriptscriptstyle {#2}{#3}}%
}
\def\MT_phantom_word:NNnn #1#2#3#4{%
\@begin@tempboxa\hbox{$\m@th#2#4$}%
\setlength\@tempdima{\widthof{$\m@th#2#3$}}%
\hbox{\hb@xt@\@tempdima{\csname bm@#1\endcsname}}%
\@end@tempboxa}
\MHInternalSyntaxOff

\usepackage{amsmath,amsfonts,amssymb}
\usepackage{graphicx}
\usepackage[colorlinks=true, allcolors=blue]{hyperref}

\title{Nonlinear  modeling and  preliminary stabilization results for a class of piezoelectric smart composite beams}

\author[a]{Ahmet \"Ozkan \"Ozer }
\affil[a]{Department of Mathematics, Western Kentucky University, Bowling Green, KY, 42101,  USA}

\authorinfo{Ahmet \"{O}zkan \"{O}zer is a faculty member in the Department of Mathematics of Western Kentucky University, Bowling Green, KY, 42101,  USA
        {\tt\small ozkan.ozer@wku.edu}}

\begin{document}
\maketitle

\begin{abstract}
Existing smart composite piezoelectric beam  models in the literature mostly ignore the electro-magnetic interactions and adopt the linear elasticity theory. However,  these interactions  substantially change the controllability and stabilizability at the high frequencies, and linear models fail to represent and predict the governing dynamics since mechanical nonlinearities are pronounced  in certain applications such as energy harvesting.
In this paper, first, a consistent variational approach is used by considering nonlinear elasticity theory to derive equations of motion for a single-layer  piezoelectric beam with and without the electromagnetic interactions (fully dynamic and  electrostatic). This modeling strategy  is extended for the three-layer piezoelectric smart composites by adopting the two widely-accepted sandwich beam theories. For both single-layer and three-layer models, the resulting infinite dimensional equations of motion can be formulated in the state-space form. It is observed that the fully dynamic nonlinear models are unbounded boundary control systems (same in linear theory) ${\bf \dot y}(t)=(\mathcal A +\mathcal N) {\bf y} (t) + \mathcal B u(t)$, the electrostatic nonlinear models are unbounded bilinear control systems ${\bf \dot y}(t)=(\mathcal A +\mathcal N){\bf y} (t) + (\mathcal B_1+ \mathcal B_2 {\bf y}) u(t)$ in sharp contrast to the linear theory. Finally, we propose $\mathcal B^*-$type feedback controllers to stabilize the single piezoelectric  beam models. The  filtered semi-discrete Finite Difference approximations is adopted to illustrate the findings.
\end{abstract}


\section{Introduction}
\begin{figure}
    \centering
    \begin{subfigure}[t]{0.45\textwidth}
        \centering
        \includegraphics[width=\linewidth]{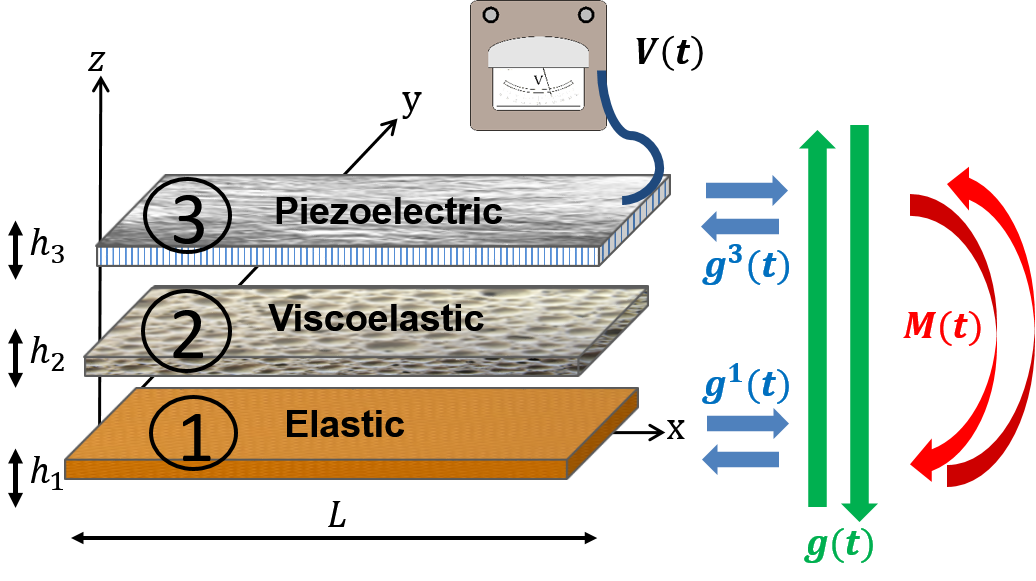}
        \caption{\scriptsize A piezoelectric smart composite of length $L$  with thicknesses $h_1, h_2, h_3$ for its layers {\normalsize \textcircled{1},  \textcircled{2}, \textcircled{3},} respectively. Voltage control $\bf V(t)$  and strain controllers $\bf g^1(t), g^3(t)$  (in $x-$direction) control stretching motion of the layers  {\normalsize \textcircled{1} and  \textcircled{3}}. The whole device bends as  a whole, and the bending motions are controlled by  shear $\bf g(t)$ and moment  $\bf M(t)$ controllers. As these controllers are actuated, vibrational modes on the device can be suppressed in a fraction of a second. } \label{controls2}
    \end{subfigure}
    \hfill
    \begin{subfigure}[t]{0.45\textwidth}
        \centering
        \includegraphics[width=\linewidth]{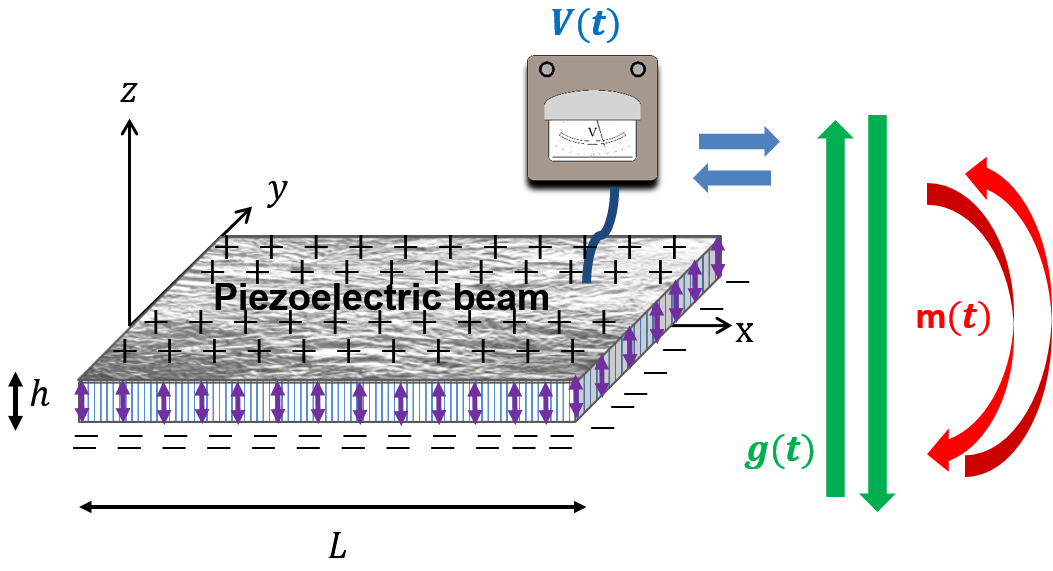}
        \caption{\scriptsize  A single piezoelectric beam of length $L$  with thickness $h.$  Voltage control $\bf V(t)$ at the circuit and strain controller $\bf g(t)$  (in $x-$direction) control stretching motion. Since the nonlinear theory is used, the bending motions are controlled by  shear $\bf g(t)$ and moment  $\bf M(t)$ controllers.} \label{single}
    \end{subfigure}
    \end{figure}

A perfectly-bonded three-layer piezoelectric smart composite beam  is an elastic structure consisting of perfectly bonded  piezoelectric and elastic layers  {\large \textcircled{1}}  and {\large \textcircled{3},} and a viscoelastic layer {\large \textcircled{2}}, see Fig. \ref{controls2}. The layer {\large \textcircled{2}} may be foam-based and it has the ability to passively suppress the vibrations on the device \cite{Ozkan2}. These devices have become more and more promising in industrial applications in aeronautic, civil, defense, biomedical, and space structures  \cite{Dag1,Smith} due to their small size
and high power density. Applications include cardiac pacemakers \cite{Dag}, course-changing bullets \cite{Exacto}, structural health monitoring \cite{health},  nano-positioners \cite{Gu}, ultrasound imagers \cite{Sci-D}, ultrasonic welders and cleaners \cite{Ultra}, energy harvesting \cite{E-Inman} due to the excellent advantages of the fast response time, large
mechanical force, and extremely fine resolution \cite{Gu,Ru}.  Accurate modeling and controlling of these structures are critical to achieve objectives for high-precision motion \cite{Cao-Chen}.

As the linear mechanical effects are considered, certain mechanical and electro-magnetic assumptions are required for each layer. The middle layer is not resistant to shear, and outer layers are relatively more stiff. Therefore,  Mindlin-Timoshenko and Euler-Bernoulli-type linear elasticity theories \cite{Lagnese-Lions} are good fits for the middle and outer layers, respectively. For the electro-magnetic interactions on the piezoelectric layer,  electrostatic assumption is almost standard through existing models \cite{Baz,Lam,Shen}.   For the mechanical interaction of layers  {\large \textcircled{1}}, {\large \textcircled{2}} and {\large \textcircled{3},}   the existing piezoelectric composite beam models  use  Mead-Marcus \cite{Mead} or  Rao-Nakra \cite{Rao} sandwich beam theories. For instance,  a Mead-Marcus type model is obtained by neglecting the rotational inertia terms for the longitudinal dynamics and rotational inertia for the bending dynamics \cite{Baz}. On the other hand, the model obtained in \cite{Lam,Shen} through a consistent variational approach is more like a Rao-Nakra-type.  All of these models  reduce to the classical counterparts (sandwich laminates)  once the piezoelectric strain is taken to be zero \cite{Hansen3}.
Recently, linear fully dynamic and electrostatic models are obtained through a consistent variational approach \cite{O-M}-\cite{Ozkan7}. The electrostatic models are simply derived by freeing the magnetic effects  in the fully dynamic models.  Keeping magnetic effects in the equations  provide not only the natural control and observation operators but also exploring how the electrical and mechanical equations are coupled rigorously. In fact, the $\mc B^*-$observation is the total induced current at the electrodes \cite{O-M1,Ozkan3}. It  is much  easier and physical in terms of practical applications since measuring the total induced current at the electrodes is easier than measuring displacements or the velocity of the composite at one end of the beam in the electrostatic case, i.e. see \cite{Baz,Ozkan2,Ronkanen,Smith}. As well, there is no need of an external displacement or velocity sensor \cite{Fur}.

As the nonlinear mechanical effects are considered, the electro-magnetic behavior for a single piezoelectric layer, see Fig. \ref{single}, is still modeled by the electrostatic assumption in the literature  \cite{Baz,Kugi}.  In deep contrast to the linear elasticity theory, the voltage control actuating the piezoelectric layer can  control both the stretching and bending. Moreover, the equations of motion form an unbounded bilinear system. This is first observed in \cite{Kugi}. Even though it is known that the electrostatic theory is sufficient for many applications in piezoelectric acoustic wave devices, there are situations in which full electromagnetic coupling needs to be considered \cite{Yang}.
To our knowledge, no work is reported yet on the fully-dynamic nonlinear vibrations of piezoelectric beams other than the port-Hamiltonian framework \cite{Voss} where it is shown that the quasi-static model lacks of asymptotic stability.  For the nonlinear piezoelectric sandwich beams, an electrostatic model is obtained by adopting the  Mead-Marcus sandwich theory. To our knowledge, the literature is scarce of rigorous nonlinear  electrostatic and fully dynamic piezoelectric composite beam models and related stabilization results.

% The electromagnetic modeling of the piezoelectric layer  traditionally uses  the electrostatic assumption which rules out  all dynamic electrical  and magnetic effects due to Maxwell's equations, i.e. see  \cite{Smith}.
  %This means the electromagnetic equations are static, and the electric field and the magnetic field are not dynamically coupled so the wave behavior of electromagnetic fields is not seen.

%The goal of the research is to (i) obtain nonlinear models for a single piezoelectric beam with/without magnetic effects by using the large displacement assumption of Euler Bernoulli and Mindlin-Timoshenko theories \cite{Lagnese-Lions}, (ii) obtain nonlinear piezoelectric smart composite beam models by using the Mead-Marcus and Rao-Nakra sandwich beam theory assumptions with/without magnetic effects, (iii)  show that the electrostatic models for a nonlinear single piezoelectric beam and the Rao-Nakra type piezoelectric smart composite beam  are ``unbounded" bilinear (affine) control systems with the voltage input, (iv)  provide a foundation so that it is possible to analyze the well-posedness, controllability and stabilization problems.

The paper is organized as the following. In Section \ref{piezo}, both dynamic and electrostatic nonlinear models for a single cantilevered piezoelectric beam are derived through a variational approach by using the Euler Bernoulli and Mindlin-Timoshenko theories. In Section \ref{S3}, both dynamic and electrostatic  nonlinear models for a three-layer piezoelectric smart cantilevered beam are derived  by using  Mead-Marcus and Rao-Nakra beam theories. In both sections, the electrostatic models are rigorously obtained by freeing the magnetic effects in the differential equations. All electrostatic models are ``unbounded" bilinear (affine) control systems with the voltage input. These models form  a solid foundation so that it is possible to analyze the well-posedness, stabilization, and approximation problems. In Section \ref{Num}, we propose $\mc B^*-$type (nonlinear) stabilizing controllers for the closed-loop systems of electrostatic piezoelectric beam models. Finally, the  filtered semi-discrete Finite Difference approximation is adopted to illustrate the findings. All of the  models obtained in the paper can be also used for the dual problem ``energy harvesting" by simply coupling the equations to a resistive circuit equation \cite{E-Inman,Shubov}, see Remark \ref{har}.

\section{Fully dynamic single layer piezoelectric beam }
\label{piezo}
Let $x, y$ and $z$ be the longitudinal and transverse directions, respectively, and the piezoelectric beam occupy the region $\Omega=[0,L]\times [-b,b]\times [-\frac{h}{2}, \frac{h}{2}]$  where $h << L ,$ see Fig. \ref{single}. A very widely-used linear constitutive relationship  for piezoelectric beams  is
$\left( \begin{array}{l}
 T \\
 D \\
 \end{array} \right)=
\left[ {\begin{array}{*{20}c}
   \alpha & -\gamma^{\text{T}}  \\
   \gamma & \varepsilon  \\
\end{array}} \right]\left( \begin{array}{l}
 S \\
 E \\
 \end{array} \right)$
where $T=(T_{11}, T_{22}, T_{33}, T_{23}, T_{13}, T_{12})^{\text T},$ $S=(S_{11}, S_{22}, S_{33}, S_{23}, S_{13}, S_{12})^{\text T}$ $D=(D_1, D_2, D_3)^{\text T},$ and  $E=(E_1, E_2, E_3)^{\text{T}}$ are stress, strain, electric displacement, and electric field vectors, respectively. Moreover, the matrices $[\alpha]_{6\times 6}, [\gamma]_{3\times 6},$ $ [\varepsilon]_{3\times 3}$ are the matrices with elastic, electro-mechanic, and dielectric constant entries (for more details the reader can refer to \cite{Tiersten}).   We assume the transverse isotropy and polarization in $x_3-$direction. Since $h << L$,  assume that all forces acting in the $x_2$ direction are zero. We also assume that  $ T_{33}$ is zero, and therefore
$$T=(T_{11}, T_{13})^{\text T}, S=(S_{11}, S_{13})^{\text T}, D=(D_1, D_3)^{\text T}, E=(E_1, E_3)^{\text{T}}.$$
Let $v=v(x),$ $w=w(x)$ and $\psi=\psi(x)$ denote the longitudinal displacement of the center line, transverse    displacement of the beam, and the rotation, respectively. We assume the Euler-Bernoulli (E-B) and Mindlin-Timoshenko (M-T) large-displacement assumptions related nonzero components of strains $$S_{11}=\frac{\partial U_1}{\partial x} + \frac{1}{2}\left(\frac{\partial U_3}{\partial x}\right)^2,~~ S_{13}=\frac{\partial U_1}{\partial z}+ \frac{\partial U_3}{\partial x} + \frac{1}{2}\frac{\partial U_3}{\partial x}  \frac{\partial U_3}{\partial z}$$ where it is assumed that $\frac{\partial U_1}{\partial x}$ is small relative to the $\frac{\partial U_3}{\partial x}.$  Defining the physical variables
    \begin{eqnarray}
    \begin{array}{cc}
    \label{coef} \gamma_1=\gamma_{15}, ~~\gamma_3=\gamma_{31},~~\beta_{1}=\frac{1}{\varepsilon_{11}}, ~~\beta_3=\frac{1}{\ep3}~~ \alpha_1=\alpha_{11} + \gamma_3^2 \beta_3, ~~\alpha_3=\alpha_{55}+ \gamma_1^2 \beta_1.
 \end{array}
 \end{eqnarray} the displacement fields for (E-B) and (M-T) assumptions, and the corresponding constitutive relationship between $(T,E)$ and $(S,D)$ are obtained  in Table \ref{table:displ}.
 \bgroup
\def\arraystretch{1.1}
\begin{table}[h]
\centering
\begin{tabular}{|p{5cm}|p{5cm}|} \hline
Euler-Bernoulli (E-B) & Mindlin-Timoshenko (M-T) \\ \hline\hline
$U_1=v(x)-zw_x(x)$ & $U_1=v(x)+z\psi(x)$  \\ \hline
$U_3=w(x)$ & $U_3=w(x)$  \\ \hline\hline
$S_{11}=v_x-zw_{xx}+\frac{1}{2}w_x^2$ & $S_{11}=v_x+z\psi_x +\frac{1}{2}w_x^2$  \\ \hline
$S_{13}=0$ & $S_{13}=\frac{1}{2} (w_x+ \psi)$  \\ \hline\hline
$T_{11}=\alpha_1 S_{11}-\gamma_3\beta_3 D_3$ & $ T_{11}= \alpha_1 S_{11}-\gamma_3\beta_3 D_3$  \\ \hline
$T_{13}=-\gamma_{1}\beta_{1} D_1$ & $ T_{13}= \alpha_3\ S_{13}-\gamma_1 \beta_1 D_1 $  \\ \hline
$E_1=\beta_{1}D_1$ & $E_1=-\gamma_1 \beta_1 S_{13}+\b1 D_1$  \\ \hline
$E_3=-\gamma_3\beta_3 S_{11}+\beta_3 D_3$ & $E_3=-\gamma_3\beta_3 S_{11}+\bp3 D_3$  \\ \hline
\end{tabular}
\caption{\small Nonlinear constitutive relationships for (E-B) and (M-T) assumptions.}
\label{table:displ}
\end{table}

 The following Lagrangian is used for voltage-driven piezoelectric beams so that the applied voltage appears in the work term:
    \begin{eqnarray}\label{tildeL}  {\mb L}= \int_0^T \left[\mb{K}-(\mb{P}+\mb{E})+\mb B +\mb{W}\right]~dt\end{eqnarray}
 where $\mb K, \mb P, \mb E$ and  $\mb B$  denote kinetic, potential, electrical, magnetic energies of the beam, respectively, and $\mb W$ is the work done by the external forces. Moreover, $\mb P+\mb E$ is the total stored energy of the beam \cite{O-M1}. % The Lagrangian ${\mb L}$ is a function of independent variables $(S, D).$

The magnetic energy $\mb B$ (as a function of $D$) is added to the Lagrangian ${\mb L}$ through Maxwell's equations. Let $B$ denote magnetic field vector
In this paper, we first use the dynamic approach for the modeling of a piezoelectric beam. Maxwell's equations including the effects of $\textrm{B}$ are
$ \nabla\cdot \textrm{B}=0, ~~~\dot {\textrm{B}}= -\nabla\times E,$ and $ \dot D= \frac{1}{\mu}(\nabla\times \textrm{B})$ where $\mu$ is the magnetic permeability. It follows from the last equation that  $\frac{ d \textrm{B}_2}{dx}=-\mu \dot D_3,$ and so $\textrm{B}_2=-\mu\int_0^x \dot D_3(\xi, x_3, t) ~d\xi.$
%The magnetic energy, which can be regarded as the   ``electrical kinetic energy'', is
%\begin{eqnarray*}
%\mb B= \frac{1}{2\mu}\int_{\Omega} \|\textrm{B}\|^2~dX =  \frac{\mu}{2}\int_{\Omega} \left[\int_0^x \dot D_3(\xi, x_3, t) d\xi\right]^2 ~dX.
%\end{eqnarray*}
Define $p=\int_0^x  D_3(\xi,  t) ~d\xi$ to be the total charge of the beam at $x\in (0,L]$ so that $p_x=D_3.$ By using Table \ref{table:displ} (with the assumption $D_1=0$), magnetic (electrical kinetic), kinetic, and stored energies of the beam are given respectively by
\begin{eqnarray}
\nonumber &&\mb{B}= \frac{\mu h}{2}\int_{0}^L \dot p^2~dx,~ \text{(E-B)}, \text{(M-T)}\\
  \nonumber &&\mb{K}=\frac{\rho}{2} \int_\Omega \left(\dot U^2+ \dot W^2\right)~dX=\frac{\rho h}{2} \left\{
  \begin{array}{ll}
    \int_0^L \left(\dot v^2+\frac{ h^2}{12} \dot w_x^2 + \dot w^2\right)~dx, & \text{(E-B)}\\
\int_0^L \left(\dot v^2+\frac{ h^2}{12} \dot \psi^2 + \dot w^2\right)~dx,& \text{(M-T)}\\
  \end{array}
\right.\\
\nonumber && \mb P+\mb E =\frac{1}{2}\int_\Omega \left(T_{11}S_{11} + T_{13}S_{13} + D_3 E_3\right) ~dX =\left\{
  \begin{array}{ll}
     \frac{h}{2}\int_0^L \left[ \alpha_1  \left( \left(v_x + \frac{1}{2}w_x^2\right)^2+ \frac{ h^2}{12}w_{xx}^2 \right) \right., &  \\
   \quad \quad \left.-2{\gamma_3\beta_3}  \left(v_x+ \frac{1}{2}w_x^2\right) p' + \beta_{3} p_x^2\right]~dx,  \text{(E-B)} &\\
    \frac{h}{2} \int_0^L \left[ \alpha_1  \left( \left(v_x + \frac{1}{2}w_x^2\right)^2+ \frac{ h^2}{12}\psi_x^2 \right)+ \beta_{3} p_x^2 \right.&\\
\left.\quad +\alpha_3(w_x+\psi)^2-2{\gamma_3\beta_3}  \left(v_x+ \frac{1}{2}w_x^2\right) p_x \right]~dx,~ \text{(M-T)} &
%\quad\quad \left.+ \beta_{3} p_x^2\right]~dx,\quad\quad \text{(M-T)}
  \end{array}
\right.
\end{eqnarray}
where $\rho,\mu$ denote the volume density and magnetic permeability of the beam, respectively. Assume that the beam is subject to a distribution of boundary forces $(\tilde g^1, \tilde g)$  along its edge $x=L.$  Now define
\begin{eqnarray}
 \nonumber   \left.
\begin{array}{ll}
% f^1(x,t)=\int_{0}^{z_1} \tilde f^1(x,z,t)~dz, ~f^3(x,t)=\int_{z_2}^{h} \tilde f^3(x,z,t)~dz,  & \\
\hspace{-0.1in} g^1(x,t )=\int_{-h/2}^{h/2} \tilde g^1(x,z,t)~dz,  g(x,t )=\int_{-h/2}^{h/2} \tilde g(x,z,t)~dz,
 ~~m(x,t)=\int_{-h/2}^{h/2} z~ \tilde g^1(x,z,t)~dz
 \end{array} \right.
\end{eqnarray}
to be the external force resultants defined as in \cite{Lagnese-Lions}. For our model, it is appropriate to assume that  $\tilde g^1$ is independent of $z,$ Let $V(t)$ be the voltage applied at the electrodes of the piezoelectric layer.
The total work  done by the external forces is ${\mb W}= {\mb W}_{e} + {\mb W}_{m}$ where ${\mb W}_e=\int_0^L \left[- p_x V(t) \right]~dx$ denotes the electrical force generated by applying voltage $V(t)$ to the electrodes of the piezoelectric layer, and ${\mb W}_m$ denotes the external mechanical forces  and defined  by
\begin{eqnarray}
\nonumber {\mb{W}}_{m} =\left\{
  \begin{array}{ll}
    g^1 v^1 (L)    +  g w(L)- m w_x(L), & \text{(E-B)}\\
{  g^1 v^1 (L)    +  g w(L)- m \psi(L)}.& \text{(M-T)}\\
  \end{array}
\right.
%\label{work-done} \mb{W}= - \int_{S} D_3 ~\bar\phi ~d\Gamma=-\int_0^L  p' V(t))~dx,\quad \text{(E-B)}, \text{(M-T)}.
\end{eqnarray}
 %where $N=\int_{-h/2}^{h/2} T_{11} dz =c_{11}^D h\left(v_x + \frac{1}{2}w_x^2\right)$ is the stress resultant,
The application of Hamilton's principle, using clamped-free boundary conditions and  setting the variation of admissible displacements $\{v,w,p\}$ for (E-B) and $\{v,w, \psi, p\}$ for  (M-T)  of ${\mb L}$ to zero, respectively yields
\begin{eqnarray}
  \label{homo-volll-EB}
  &&(E-B)\left\{
\begin{array}{ll}
  \rho h \ddot v- \alpha_1 h \left(v_x+ \frac{1}{2}w_x^2\right)_x+{\gamma_3\beta_3}  h   p_{xx} = 0,   & \\
  \mu h \ddot p   -\beta_{3} h p_{xx} + {\gamma_3\beta_3}h \left(v_x+ \frac{1}{2}w_x^2\right)_x= 0, &\\
   \rho h \ddot w -\frac{\rho h^3}{12}\ddot w_{xx}+\frac{\alpha_1 h^3}{12}w_{xxxx}+\left(\left(-\alpha_1 h \left(v_x + \frac{1}{2}w_x^2\right) + {\gamma_3\beta_3}h p_x\right) w_x\right)_x=0, &\\
     &\\
   \left[ v,p,w,w_x\right](0)=0, ~ \left[\frac{\alpha_1h^3}{12}w_{xx}\right](L)=-m(t),\quad  \left[ \alpha_1 h \left(v_x+ \frac{1}{2}w_x^2\right) - {\gamma_3\beta_3}h p_x\right](L) =g^1(t),&\\
  \left[  \beta_{3} h p_x-{\gamma_3\beta_3}h  \left(v_x+ \frac{1}{2}w_x^2\right)\right](L)=-V(t), &\\
 \left[\frac{\rho h^3}{12}\ddot w_x-\frac{\alpha_1 h^3}{12}w_{xxx}+\left(\alpha_1 h\left(v_x + \frac{1}{2}w_x^2\right)-{\gamma_3\beta_3}h p_x\right)w_x\right](L) =g(t),&\\
 \label{initial-EB}  (v, w, p, \dot v, \dot w, \dot p)(x,0)=(v_0, w_0, p_0, v_1, w_1, p_1),&
\end{array}\right.\\
  \label{homo-volll-MT}
&& (M-T) \left\{
\begin{array}{ll}
  \rho h \ddot v- \alpha_1 h \left(v_x + \frac{1}{2}w_x^2\right)_x+{\gamma_3\beta_3}  h   p_{xx} = 0,   & \\
   \frac{\rho h^3}{12} \ddot \psi -\frac{\alpha_1 h^3}{12}\psi_{xx}+\alpha_3 h\left(w_x+ \psi\right)=0, &\\
  \rho h \ddot w - \alpha_1 h\left(\left(v_x + \frac{1}{2}w_x^2\right) w_x\right)_x+{\gamma_3\beta_3}h (w_x p_x)_x-\alpha_3 h\left(w_x + \psi\right)_x=0, &\\
  \mu h \ddot p   -\beta_{3} h p_{xx} + {\gamma_3\beta_3}h \left(v_x+ \frac{1}{2}w_x^2\right)_x= 0, &\\
  &\\
    \left[v,p,w,\psi\right](0)=0,\quad \frac{\alpha_1 h^3}{12}\psi_{x}(L)=m(t),  \quad  \left[ \alpha_1 h \left(v_x + \frac{1}{2}w_x^2\right) - {\gamma_3\beta_3}h p_x\right]_{x=L} =g^1(t),&\\
    \left[\beta_{3} h p_x-{\gamma_3\beta_3}h  \left(v_x + \frac{1}{2}w_x^2\right)\right]_{x=L}=-V(t), &\\
\left[ \alpha_1 h \left(v_x+ \frac{1}{2}w_x^2\right)w_x- {\gamma_3\beta_3}h w_x p_x + \alpha_3 h (w_x+\psi)\right](L)=g(t), &\\
 (v, w, \psi, p, \dot v, \dot w, \dot\psi, \dot p)(x,0)= (v_0, w_0, \psi_0, p_0, v_1, w_1, \psi_1, p_1).&
\end{array}\right.
\end{eqnarray}
Note that these nonlinear equations are novel and  they were never studied before. If one considers the linearization of (\ref{homo-volll-EB}) and (\ref{homo-volll-MT}) around  $v=w=p=\dot v=\dot w=\dot p=0$ for (E-B) and  $v=\psi=w= p=\dot v=\dot \psi=\dot w=\dot p=0$ for (M-T),  or if we simply consider the linear stress-strain relationship in the beginning, we still obtain a  boundary control system \cite{O-M1}}
\begin{eqnarray}
 && (E-B)\left\{
\begin{array}{ll}
  \rho h \ddot v- \alpha_1 h v_{xx}+{\gamma_3\beta_3}  h   p_{xx} = 0,   & \\
  \mu h \ddot p   -\beta_{3} h p_{xx} + {\gamma_3\beta_3}h v_{xx}= 0, &\\
   \rho h \ddot w -\frac{\rho h^3}{12}\ddot w_{xx}+\frac{\alpha_1 h^3}{12}w_{xxxx}=0, &\\
     &\\
   \left[ v,p,w,w_x\right](0)=0, \quad  \left[\frac{\alpha_1h^3}{12}w_{xx}\right](L)=-m(t),\quad \left[ \alpha_1 h v_x- {\gamma_3\beta_3}h p_x\right](L) =g^1(t),&\\
  \left[  \beta_{3} h p_x-{\gamma_3\beta_3}h  v_x\right](L)=-V(t), \quad  \left[\frac{\rho h^3}{12}\ddot w_x-\frac{\alpha_1 h^3}{12}w_{xxx}\right](L) =g(t),&\\
 (v, w, p, \dot v, \dot w, \dot p)(x,0)=(v_0, w_0, p_0, v_1, w_1, p_1),&
\end{array}\right.\\
 \label{initial-EB}
 && (M-T) \left\{
\begin{array}{ll}
  \rho h \ddot v- \alpha_1 h v_{xx}+{\gamma_3\beta_3}  h   p_{xx} = 0,   & \\
   \frac{\rho h^3}{12} \ddot \psi -\frac{\alpha_1 h^3}{12}\psi_{xx}+\alpha_3 h\left(w_x+ \psi\right)=0, &\\
  \rho h \ddot w -\alpha_3 h\left(w_x + \psi\right)_x=0, &\\
  \mu h \ddot p   -\beta_{3} h p_{xx} + {\gamma_3\beta_3}h \left(v_x+ \frac{1}{2}w_x^2\right)_x= 0, &\\
  &\\
    \left[v,p,w,\psi\right](0)=0,\quad     \frac{\alpha_1 h^3}{12}\psi_{x}(L)=m(t),  \quad    \left[ \alpha_1 h v_x - {\gamma_3\beta_3}h p_x\right]_{x=L} =g^1(t),&\\
    \left[\beta_{3} h p_x-{\gamma_3\beta_3}h  v_x\right]_{x=L}=-V(t), \quad  \left[  \alpha_3 h (w_x+\psi)\right](L)=g(t), &\\
 (v, w, \psi, p, \dot v, \dot w, \dot\psi, \dot p)(x,0)= (v_0, w_0, \psi_0, p_0, v_1, w_1, \psi_1, p_1).&
\end{array}\right.
\end{eqnarray}
  The bending equation in (E-B) and the bending and rotation angle equations in (M-T)  are completely decoupled from the rest. If one considers the voltage control $V(t)$ only, the bending and shear motions can never be controlled. Only certain longitudinal motions are  controlled \cite{O-M,O-M1} by measuring the  total induced current accumulated at the electrodes \cite{Ozkan3}. In practical applications, these uncontrolled longitudinal motions corresponding to the high-frequency solutions.
\subsection{Electrostatic models} By the electrostatic assumption, all magnetic effects are discarded  (i.e., take $\ddot p\equiv 0$) and  $g^1(t)\equiv 0.$ Therefore  (\ref{homo-volll-EB}) and (\ref{homo-volll-MT}) respectively reduce to
\begin{eqnarray}
  \label{homo-NoM-EB}
 (E-B) \left\{
\begin{array}{ll}
  \rho h \ddot v- \alpha_{11} h \left(v_x + \frac{1}{2}w_x^2\right)_x= 0   & \\
  \rho h \ddot w -\frac{\rho h^3}{12}\ddot w_{xx}+\frac{\alpha_{1}h^3}{12}w_{xxxx}- \left[\left(\alpha_{11} h\left(v_x + \frac{1}{2}w_x^2\right) + \gamma_3 H_L V(t) \right)w_x\right]_x=0, &\\
  &\\
   \left[v,w,w_x\right](0)=0, \quad \frac{\alpha_1h^3}{12}w_{xx}(L)=-m(t),&\\
    \left[ \alpha_{11}h \left(v_x+ \frac{1}{2}w_x^2\right)\right](L) =-\gamma_3 V(t), ~~ \left[\frac{\rho h^3}{12}\ddot w_x-\frac{\alpha_1 h^3}{12}w_{xxx}\right](L) =g(t),&\\
  (v, w, \dot v, \dot w)(x,0)=(v_0, w_0, v_1, w_1).&
\end{array}\right.&&\\
 \label{homo-NoM-MT}
(M-T)   \left\{
\begin{array}{ll}
  \rho h \ddot v- \alpha_{11} h \left(v_x + \frac{1}{2}w_x^2\right)_x= 0,   & \\
   \frac{\rho h^3}{12} \ddot \psi -\frac{\alpha_1 h^3}{12}\psi_{xx}+\alpha_3 h\left(w_x + \psi\right)=0, &\\
  \rho h \ddot w  -\alpha_3 h\left(w_x + \psi\right)_x  - \left[\alpha_{11} h\left(v_x + \frac{1}{2}w_x^2\right) w'+ \gamma_3 H_L V(t)w_x\right]_x=0, &\\
  &\\
  %\mu h \ddot p   -\beta_{3} h p'' + {\gamma_3\beta_3}h \left(v'+ \frac{1}{2}w'^2\right)'= 0, &\\
    \left[v,w,\psi\right](0)=0, ~~\frac{\alpha_1 h^3}{12}\psi_{x}(L)=m(t),~~ \left[\alpha_{11} h \left(v_x + \frac{1}{2}w_x^2\right)\right] (L) =-\gamma_3 V(t),&\\
    %\left[\beta_{3} h p'-{\gamma_3\beta_3}h  \left(v' + \frac{1}{2}w'^2\right)=-V(t)\right]_{x=L} &\\
  \left[ \alpha_3 h (w_x+\psi)\right](L)=g(t), &\\
 (v, w, \psi,  \dot v, \dot w, \dot\psi)(x,0)=(v_0, w_0, \psi_0, v_1, w_1, \psi_1)&
\end{array}\right.&&
\end{eqnarray}
 where we used (\ref{coef}), and $H_L=H(x)-H(x-L)$ is the characteristic function of the interval $(0,L).$ Note that the fully elastic (E-B)  and (M-T) beam models are derived and studied in \cite{Lag-L},\cite{Ashgari}, respectively by taking $V(t)\equiv 0.$ Observe that as the longitudinal inertia term $\rho \ddot v\equiv 0$ in  (\ref{homo-NoM-EB})-(\ref{homo-NoM-MT})  is neglected, the voltage control $V(t)$ does not affect the bending equation anymore. This is the difference between elastic beam models \cite{G-Li,Lag-L} and electrostatic (piezoelectric smart beam) models. The literature on which the dynamic longitudinal motions are not considered do not fit to the physics of piezoelectric beams since the bending and stretching are strongly coupled via $V(t)$.

Since piezoelectric structures are also elastic structures, damping can be added to the models above. For instance, a Kelvin-Voigt (K-V) type damping can considered  \cite{Smith}. In particular, by keeping only the linear damping terms in Table \ref{table:displ2}, the (E-B) systems (\ref{homo-volll-EB}) and (\ref{homo-NoM-EB}) can be re-written by replacing the terms $\alpha_{1} v_x$ and $\frac{\alpha_1 h^3}{12}w_{xxxx}$ by  $\alpha_{1} v_x + \tilde\alpha_{1}  \dot v_x$ and $\frac{\alpha_1 h^3}{12}w_{xxxx}+ \frac{\tilde \alpha_1 h^3}{12}\dot w_{xxxx},$ respectively, where $\tilde \alpha_{1}$ is the damping coefficient.
\bgroup
\def\arraystretch{1.2}
\begin{table}[h]
\centering
\begin{tabular}{|p{5cm}|p{5cm}|} \hline
(E-B) with (K-V) damping &(M-T) with (K-V) damping\\ \hline\hline
$T_{11}=\alpha_1 S_{11}+\tilde \alpha_1 \dot S_{11}-\gamma_3\beta_3 D_3$ & $ T_{11}= \alpha_1 S_{11}+\tilde\alpha_1 \dot S_{11}-\gamma_3\beta_3 D_3$  \\ \hline
$T_{13}=-\gamma_{1}\beta_{1} D_1$ & $ T_{13}= \alpha_3\ S_{13}+\tilde \alpha_3\dot S_{13}-\gamma_1 \beta_1 D_1 $  \\ \hline
%$E_1=\beta_{1}D_1$ & $E_1=-\gamma_1 \beta_1 S_{13}+\b1 D_1$  \\ \hline
%$E_3=-\gamma_3\beta_3 S_{11}+\beta_3 D_3$ & $E_3=-\gamma_3\beta_3 S_{11}+\bp3 D_3$  \\ \hline
\end{tabular}
\caption{\small Adding Kelvin-Voigt type damping to the constitutive equations.}
\label{table:displ2}
\end{table}
Note also that the main challenge in the nonlinear electrostatic models  (\ref{homo-NoM-EB}) and  (\ref{homo-NoM-MT}) is the voltage control $V(t)$ appearing at both  the axial strain boundary condition at $x=L$ and  the $w-$equation. This makes the boundary control system not only unbounded but also bilinear and unbounded.
%$ yet it does not fit in the bilinear Euler-Bernoulli model studied in \cite{Lenhart} where the control operator is bounded.
%(i) The boundary controllability and stabilization of a nonlinear beam model (without considering the piezoelectricity) is studied in i.e. \cite{G-Li,Lag,Lag-L}, and the references therein. In fact,
To show this, consider  (\ref{homo-NoM-EB}) with  the state $x=[v, w, \dot v, \dot w]^{\rm T}$. Let $D_x=\frac{d}{dx}, $$D_x^n=\frac{d^n}{dx^n}$ for $n>1.$ Defining the operator $\mc P=(I-\frac{h^2}{12}D_x^2)^{-1}$ where $I$ is the identity operator,  an unbounded bilinear control system for (\ref{homo-NoM-EB}) is obtained as the following
\begin{eqnarray}\label{s-s}\dot {\bf y}=(\mc A +\mc N){\bf y} +  \left(\mc B_1+\mc B_2 {\bf y}\right) u(t)\quad {\rm where}\end{eqnarray}
$$\mc A =
  \left(
             \begin{array}{c;{2pt/2pt}c}
                \begin{matrix}\frac{\alpha_{11}}{\rho}D_x^2 & 0  \\
               0&  \frac{\alpha_{11}}{\rho}\mc P\left(-\frac{h^2}{12}D_x^4+D_x^2\right)\end{matrix} &  0_{2\times 2}\\ \hdashline
               \block(1,1){0_{2\times 2}}   & \block(1,1){I_{2\times 2}}
             \end{array}
           \right),\quad \mc B_1= \left(
             \begin{array}{ccc}
             0 & 0 &0\\
             0 & 0 &0\\
           - \frac{\gamma_3D_x H_L}{\rho h} & 0 & 0\\
            0& \frac{-\mc PD_x\delta_L}{\rho h}& \frac{-\mc P\delta_L}{\rho h}
             \end{array}
           \right),$$ $\mc B_2=  \left(
             \begin{array}{ccc}
             0 & 0 &0\\
             0 & 0 &0\\
            0& 0 & 0\\
            \frac{\gamma_3}{\rho h} \mc P D_x(H_LD_x) & 0 & 0
             \end{array}
           \right),$
$u(t)=\left[
                                                                                                                     \begin{array}{c}
                                                                                                                       V(t) \\
                                                                                                                      m(t) \\
                                                                                                                      g(t)
                                                                                                                     \end{array}
                                                                                                                   \right],$ and $\mc N {\bf y}=  \frac{\alpha_{11}}{\rho}\left(
                                                                                                                     \begin{array}{c}
                                                                                                                       0_{2\times 1} \\
                                                                                                                      D_x \left(\frac{1}{2}\left(D_x(\cdot)\right)^2\right) \\
                                                                                                                        \mc PD_x\left(\frac{1}{2}\left(D_x(\cdot)\right)^3\right) \\
                                                                                                                     \end{array}
                                                                                                                   \right)$ is a nonlinear operator, and $\delta_L=\delta(x)-\delta(x-L)$. Define the natural energy space as $$\mathrm H=H^1_L(0,L)\times H^2_L(0,L)\times L^2(0,L)\times H^1_L (0,L)$$ where
 %                                                                                                                   \begin{eqnarray}
 %\nonumber
  %\begin{array}{ll}
 %H^1_L(0,L)=\{\psi\in H^1(0,L): \psi(0)=0\}, &\\
 %  H^2_L(0,L)=\{\psi\in H^2(0,L): \psi(0)=\psi_x(0)=0\},   &
  % \end{array}
%\end{eqnarray}
$$H^1_L(0,L)=\{z\in H^1(0,L): z(0)=0\},\quad H^2_L(0,L)=\{z \in H^2(0,L): z(0)=z_x(0)=0\}.$$ It can be shown that $\mc A$ is an infinitesimal generator of a unitary semigroup on $\mathrm H,$ $\mc N:\mathrm H\to \mathrm H$ is locally Lipschitz. As well,   $\mc B_1: \mathbb{C}\to \mathrm H, $  $\mc B_2:\mathbb{C}\to \mathrm H $ are unbounded and bounded operators, respectively. The local controllability or the uniform (or asymptotic) stabilization of systems like (\ref{s-s}) is still an open problem. A similar framework is recently studied in \cite{Ay,Bea,Ber}. In Section \ref{Num}, we take a stab of this hard problem in a numerical sense.

Unlike the electro-static models, the fully dynamic models (\ref{homo-volll-EB}),(\ref{homo-volll-MT}) with the choice of states ${\bf y}=[v, p, w, \dot v, \dot p,\dot w]^{\rm T}$ can be written as
$$\dot {\bf y}=(\mc A +\mc N){\bf y} + \mc B u(t)$$
 where $B$ is a boundary control operator with the input $u(t)=[g^1(t), V(t), m(t), g(t)]^{\rm T}$. Then observation $\mc B^*{\bf y}$ for only the voltage input corresponds to the total induced current $\left(\dot p(L)=\int_0^L \dot D_3 dx\right)$ accumulated at the electrodes \cite{Ozkan3}. This is completely electro-magnetic, and is more practical than measuring tip velocity at the tip of the beam. Hence, it is more promising in energy harvesting applications \cite{Dietl}.

If we linearize (\ref{homo-NoM-EB}) along the equilibrium $v=w=\dot v=\dot w=0$ for (E-B) and  $v=\psi=w=\dot v=\dot \psi=\dot w=0$ for (M-T), or if we simply consider the linear stress-strain relationship in the beginning, we obtain a ``degenerate" (strongly decoupled) system:
\begin{eqnarray}
\label{homo-NoM-EB-lin}&&  (E-B) \left\{
\begin{array}{ll}
  \rho h \ddot v- \alpha_{11} h v_{xx}= 0,   & \\
  \rho h \ddot w -\frac{\rho h^3}{12}\ddot w_{xx}+\frac{\alpha_{1}h^3}{12}w_{xxxx}=0, &\\
   \left[v,w,w_x\right](0)=0, \quad \alpha_{11}h v_x(L) =-\gamma_3 V(t), &\\
 \frac{\alpha_1h^3}{12}w_{xx}(L)=-m(t),\quad  \left[\frac{\rho h^3}{12}\ddot w_x-\frac{\alpha_1 h^3}{12}w_{xxx}\right](L) =g(t),&\\
  (v, w, \dot v, \dot w)(x,0)=(v_0, w_0, v_1, w_1),&
\end{array}\right.\\
 \label{homo-NoM-MT-lin}
&&(M-T)   \left\{
\begin{array}{ll}
  \rho h \ddot v- \alpha_{11} h v_{xx} = 0,   & \\
   \frac{\rho h^3}{12} \ddot \psi -\frac{\alpha_1 h^3}{12}\psi_{xx}+\alpha_3 h\left(w_x + \psi\right)=0, &\\
  \rho h \ddot w  -\alpha_3 h\left(w_x + \psi\right)_x =0, &\\
  %\mu h \ddot p   -\beta_{3} h p'' + {\gamma_3\beta_3}h \left(v'+ \frac{1}{2}w'^2\right)'= 0, &\\
    \left[v,w,\psi\right](0)=0, \quad   \alpha_{11} h v_x (L) =-\gamma_3 V(t),~~   \frac{\alpha_1 h^3}{12}\psi_{x}(L)=m(t),\quad  \left[ \alpha_3 h (w_x+\psi)\right](L)=g(t), &\\
    %\left[\beta_{3} h p'-{\gamma_3\beta_3}h  \left(v' + \frac{1}{2}w'^2\right)=-V(t)\right]_{x=L} &\\
 (v, w, \psi,  \dot v, \dot w, \dot\psi)(x,0)=(v_0, w_0, \psi_0, v_1, w_1, \psi_1).&
\end{array}\right.
\end{eqnarray}
 where the unbounded boundary control $V(t)$ only affects the stretching motion, not bending anymore.  Including all three controllers $V(t), m(t), g(t),$ both systems can be shown to be exponentially stabilizable \cite{Lag-L}.
\begin{rmk} For energy harvesting applications, the nonlinear models (\ref{homo-NoM-EB},\ref{homo-NoM-MT}) and the linear models (\ref{homo-NoM-EB-lin},\ref{homo-NoM-MT-lin}) are coupled to the extra circuit equation
\begin{eqnarray}
\frac{2\alpha_{11}\ep3 b L}{\alpha_1h} \dot V(t)+ i(t) =0
\end{eqnarray}
where $i(t)=\frac{1}{R}V(t)$ is the current generated by the piezoelectric beam, $R$ is the resistance of the attached circuit. 
\end{rmk}
\section{Fully dynamic three-layer piezoelectric smart beam}
\label{S3}
The three-layer piezoelectric smart composite beam considered in this paper consists  of a stiff layer, a complaint layer, and a piezoelectric layer, see Fig. \ref{controls2}. The composite occupies the
region $\Omega=\Omega_{xy}\times (0, h)=[0,L]\times [-b,b] \times (0,h)$ at equilibrium.
The total thickness $h$ is assumed to be small in comparison to the dimensions of $\Omega_{xy}$.
 The layers are indexed from 1 to 3 from the stiff  layer to the piezoelectric layer, respectively.

Let $0=z_0<z_1<z_2<z_3=h, $ with
 $h_i=z_i-z_{i-1}, \quad i=1,2,3.$
We use the rectangular coordinates $X=(x,y)$ to denote points in $\Omega_{xy},$ and  $(X, z)$ to denote points in $\Omega = \Omega^{\rm s} \cup \Omega^{\rm ve} \cup \Omega^{\rm p} $, where $\Omega^{\rm s}, \Omega^{\rm ve},$ and $\Omega^{\rm p}$ are the reference configurations of the stiff, viscoelastic, and piezoelectric layer, respectively, and they are given by
\begin{eqnarray*}
\begin{array}{cc}
\Omega^{\rm s}=\Omega_{xy}\times (z_0,z_1),~~ \Omega^{\rm ve}=\Omega_{xy}\times (z_1,z_2), ~~\Omega^{\rm p}=\Omega_{xy}\times (z_2,z_3).
\end{array}
\end{eqnarray*}
For $(X,z)\in \Omega,$ let $U(X,z) = ({U_1,U_2,U_3})(X, z)$ denote the displacement vector of the point
(from reference configuration). For the  beam theory,  all displacements  are assumed to be independent of $y-$coordinate, and $U_2\equiv 0.$ The transverse displacements is  $w(x,y,z)= U_3(x)= w^i(x)$ for  any $i$ and $x\in [0,L].$ Define $u^i(x,y,z)=U_1(x,0,z_i)=u^i(x)$  for $i=0,1,2,3$ and  for all $x\in (0,L).$

Define the vectors
$\vec \psi=[\psi^1,\psi^2, \psi^3]^{\rm T},$  $\vec \phi=[\phi^1,\phi^2,\phi^3]^{\rm T},$ $\vec v=[v^1, v^2, v^3]^{\rm T}$ where
\begin{eqnarray}\label{defs1}   \psi^i=\frac{u^i-u^{i-1}}{h_i}, ~\phi^i= \psi^i + w_x,  ~ v^i= \frac{u^{i-1}+u^i}{2},\quad i = 1, 2, 3,&&
\end{eqnarray}
and $\psi^{i}$   is the total rotation angles
(with negative orientation) of the deformed filament within the $i^{\rm th}$ layer in the
$x-z$ plane, $\phi^i$ is the (small angle approximation for
the) shear angles within each layer, $v^i$ is the longitudinal displacement of the center line of the $i^{\rm th}$ layer. For the middle layer, we apply
Mindlin-Timoshenko constitutive assumptions, while for the outer layers Kirchhoff displacement
assumptions. Therefore,
\begin{eqnarray}\label{defs3}\phi^1=\phi^3=0, ~  \psi^1=\psi^3=-w_x,  ~\phi^2=\psi^2+ w_x.\end{eqnarray}
Let $G_2$ be the shear modulus of the viscoelastic layer. Defining ${\hat z}^i = \frac{z^{i-1}+z^i}{2},$ and  \begin{eqnarray}
\label{coef1} \alpha^3=\alpha_1^3 + \gamma^2 \beta, ~\alpha_1^3=\alpha_{11}^3,~ \alpha^2=\alpha_{11}^2,~ \alpha^1=\alpha_{11}^1
  \end{eqnarray}
  where we use (\ref{coef}), the displacement fields, strains, and the constitutive relationships for each layer are given in Table \ref{table:displ-II}.
  \bgroup
\def\arraystretch{1.2}
\begin{table}[h]
\centering
\begin{tabular}{|p{3.5cm}|p{11cm}p{0cm}|} \hline
Layers&   Displacements, Stresses, Strains, Electric fields, Electric displacements&  \\ \hline
\multirow{3}{*}{Layer {\Large \textcircled{1}} - Elastic }&$U_1^1(x,z)=v^1(x)- (z-\hat z_1)w_x,~~U_3(x,z)=w(x)$ &  \\ \cline{2-3}
&$S_{11}=\frac{\partial v^1}{\partial x}+\frac{1}{2}(w_x)^2- (z-\hat z_1) \frac{\partial^2 w}{\partial x^2},~~S_{13}=0$ &  \\ \cline{2-3}
&$ T_{11}=\alpha_1S_{11},~~T_{13}=0$ &  \\ \cline{2-3}
\hline
\multirow{3}{*}{Layer {\Large \textcircled{2}} - Viscoelastic } &$U_1^2(x,z)=v^2(x)+ (z-\hat z_2)\psi^2(x),~~U_3(x,z)=w(x)$ &  \\ \cline{2-3}
&$ S_{11}=\frac{\partial v^2}{\partial x}- (z-\hat z_i) \frac{\partial \psi^2}{\partial x},~~S_{13}=\frac{1}{2}\phi^2$ &  \\ \cline{2-3}
&$ T_{11}=\alpha_1^2S_{11},~~T_{13}= 2G_{2} S_{13}$ &  \\ \cline{2-3}
\hline
\multirow{4}{*}{Layer {\Large \textcircled{3}} - Piezoelectric }&$U_1^3(x,z)=v^3(x)- (z-\hat z_3)w_x,~~U_3(x,z)=w(x)$ &  \\ \cline{2-3}
&$ S_{11}=\frac{\partial v^3}{\partial x}+\frac{1}{2}(w_x)^2- (z-\hat z_3) \frac{\partial^2 w}{\partial x^2},~~S_{13}=0$ &  \\ \cline{2-3}
&$ T_{11}=\alpha^3 S_{11}-\gamma\beta D_3,~~T_{13}=0$ & \\ \cline{2-3}
& $E_1=\beta_{1}D_1,~~E_3=-\gamma\beta S_{11}+\beta D_3$  &  \\ \hline
\end{tabular}
\caption{\small Nonlinear constitutive relationships for each layer.}
\label{table:displ-II}
\end{table}

Assume that the beam is subject to a distribution of boundary forces $(\tilde g^1, \tilde g^3, \tilde g)$  along its edge $x=L.$ Now define
\begin{eqnarray}
 \nonumber   \left.
\begin{array}{ll}
% f^1(x,t)=\int_{0}^{z_1} \tilde f^1(x,z,t)~dz, ~f^3(x,t)=\int_{z_2}^{h} \tilde f^3(x,z,t)~dz,  & \\
 g^i(x,t )=\int_{z_{i-1}}^{z_i} \tilde g^i(x,z,t)~dz,  ~g(x,t )=\int_{0}^{h} \tilde g(x,z,t)~dz,~~m_i=\int_{z_{i-1}}^{z_i} (z-\hat z_i) \tilde g^i(x,z,t)~dz, \quad  i=1,3

  %f(x,t )=\int_{0}^{h} \tilde f(x,z,t)~dz,
   %~g(x,t )=\int_{0}^{h} \tilde g(x,z,t)~dz, &\\
%M_i(x,t)=\int_{z_{i-1}}^{z_i} (z-\hat z_i) \tilde f^i(x,z,t) ~dz,
  \end{array} \right.
\end{eqnarray}
to be the external force resultants.  The total work  done by the external forces is ${\mb W}= {\mb W}_{e} + {\mb W}_{m}$ where ${\mb W}_e=\int_0^L \left[- p_x V(t) \right]~dx$ denotes the electrical force generated by applying voltage $V(t)$ to the electrodes of the piezoelectric layer, and ${\mb W}_m$ denotes the external mechanical forces  and defined  by
\begin{eqnarray}
\label{work-done} {\mb{W}}_m =  g^1 v^1 (L)  + g^3 v^3 (L) +  g w(L)-Mw_x(L), \quad M=m_1+m_3.
%\label{work-done} \mb{W}= - \int_{S} D_3 ~\bar\phi ~d\Gamma=-\int_0^L  p' V(t))~dx,\quad \text{(E-B)}, \text{(M-T)}.
\end{eqnarray}
 %where $N=\int_{-h/2}^{h/2} T_{11} dz =c_{11}^D h\left(v_x + \frac{1}{2}w_x^2\right)$ is the stress resultant,

 We follow the same approach in Section \ref{piezo} to include the magnetic effects. The magnetic field $\textrm{B}$ is perpendicular to the $x-z$ plane, and therefore,  $\textrm{B}_2(x)$ is the only non-zero component. Assuming  $E_1=D_1=0,$   % obtain $\textrm{B}_2=-\mu\int_0^x \dot D_3(\xi, z, t) ~d\xi.$
the  Lagrangian for the ACL beam is the same as (\ref{tildeL}) with
 %\begin{eqnarray} \nonumber \mb{L}= \int_0^T \left[\mb{K}-(\mb{P}+\mb{E})+\mb B +\mb{W}\right]~dt\end{eqnarray}
 %where
  \begin{eqnarray}
\nonumber  \begin{array}{ll}
 \hspace{-0.1in}\mb{K}= \frac{1}{2} \int_0^L \left[\sum_{i=1}^3\left(\rho_ih_i(\dot v^i)^2\right) + \left(\sum_{i=1}^3\rho_ih_i\right)(\dot w)^2 + \left(\rho_1h_1 + \rho_3h_3\right) \dot w_x^2+ \rho_2h_2 (\dot \psi^2)^2\right]~dx,&\\
 \hspace{-0.1in} \mb P+\mb E = \frac{1}{2} \int_0^L \left[\sum_{i=1,3} \left(\alpha^i h_i \left( (v_x^i + \frac{1}{2}(w_x)^2)^2 +\frac{h_i^2}{12}w_{xx}^2\right)\right)
 + \alpha^2 h_2 \left((v^2_x)^2 + \frac{h_2^2}{12}\psi_{x}^2\right)-2\gamma\beta h_3 \left( v^3_x + \frac{1}{2}(w_x)^2
 \right) p_x \right.&\\
  \quad\quad\quad \left. + \beta h_3 p_x^2 + G_2h_2 \left(\phi^2\right)^2  \right]~dx,\quad\quad&\\
  \hspace{-0.1in}\mb{B}= \frac{\mu h_3}{2}\int_{0}^L \dot p^2~dx&
% \mb{W} =\int_0^L \left[- p_x V(t) \right]~dx+  g^1 v^1 (L)  + g^3 v^3 (L) &\\
% \quad\quad +  g w(L)-Mw_x(L)
 \end{array}
\end{eqnarray}
where $\rho_i$ is the volume density of the $i^{\rm th}$ layer. By using (\ref{defs1})-(\ref{defs3}), the variables $v^2, \phi^2$  and $\psi^2$ can be written as the functions of state variables as the following
\begin{eqnarray}
\begin{array}{cc}
\nonumber v^2=\frac{1}{2}\left(v^1+v^3\right) + \frac{h_3-h_1}{4}w_x,~~\psi^2 =\frac{1}{h_2}\left(-v^1+v^3\right) + \frac{h_1 +  h_3}{2h_2}w_x,~~ \phi^2 =\frac{1}{h_2}\left(-v^1+v^3\right) + \frac{h_1 + 2h_2+ h_3}{2h_2}w_x.
\end{array}
\end{eqnarray}
 Thus, we choose $w,v^1,v^3$ as the state variables.  Let $H=\frac{h_1 + 2h_2+h_3}{2}.$ Application of Hamilton's principle, by using forced boundary conditions, i.e. clamped at $x=0,$ by setting the variation of admissible displacements $\{v^1,v^3, p, w\}$ of ${\mb L}$ in (\ref{tildeL}) to zero yields the following coupled equations of stretching in odd layers, dynamic charge in the piezoelectric layer, and the bending of the whole composite:
\begin{eqnarray}
 \label{perturbed}\left\{\hspace{-0.05in}
  \begin{array}{ll}
  \left(\rho_1h_1 + \frac{\rho_2h_2}{3}\right)\ddot v^1 + \frac{\rho_2 h_2}{6}\ddot v^3- \frac{1}{6}\alpha^2 h_2 v^3_{xx}   -\left[\left(\alpha^1 h_1 (v^1_x+\frac{1}{2}(w_x)^2\right)+ \frac{1}{3}\alpha^2 h_2 v^1_{x}\right]_x  - G_2  \phi^2 &\\
\quad\quad  -\frac{\rho_2 h_2}{12}(2h_1-h_3)\ddot w_x   + \frac{\alpha^2 h_2}{12}(2h_1-h_3) w_{xxx}= 0,&\\
  % &\\
   \left(\rho_3h_3 + \frac{\rho_2h_2}{3}\right)\ddot v^3 + \frac{\rho_2 h_2}{6}\ddot v^1- \frac{1}{6}\alpha^2 h_2 v^1_{xx} -\left[\left(\alpha^3 h_3 (v^3_x+\frac{1}{2}(w_x)^2\right)+ \frac{1}{3}\alpha^2 h_2 v^3_{x}\right]_x + G_2 \phi^2&\\
  \quad\quad  +\frac{\rho_2 h_2}{12}(2h_3-h_1)\ddot w_x     - \frac{\alpha^2 h_2}{12}(2h_3-h_1) w_{xxx} + \gamma \beta h_3 p_{xx} = 0,   & \\
   % &\\
    \mu  h_3 \ddot p   -\beta h_3   p_{xx} + \gamma \beta h_3 \left(v^3_{x}+\frac{1}{2}w_x^2\right)_x= 0, &\\
   % &\\
   (\rho_1 h_1 + \rho_2 h_2+ \rho_3 h_3) \ddot w -  H G_2   \phi^2_x   - \frac{1}{12}\left(\rho_1 h_1^3 + \rho_3 h_3^3 + \rho_2h_2(h_1^2 + h_3^2 -h_1h_3)\right) \ddot{w}_{xx} &\\
    \quad + \frac{1}{12}\left(\alpha^1 h_1^3+\alpha^3 h_3^3 + \alpha^2 h_2(h_1^2+h_3^2 -h_1h_3) \right)w_{xxxx}\hspace{-0.2in}  +\frac{\rho_2 h_2}{12} (2h_1-h_3) \ddot v^1_{x}-\frac{\alpha^2 h_2}{12}(2h_1-h_3) v^1_{xxx} &\\
   \quad  -\frac{\rho_2 h_2}{12}(2h_3-h_1) \ddot v^3_{x} +\frac{\alpha^2 h_2}{12} (2h_3-h_1) v^3_{xxx} -\sum_{i=1,3}\left[ \alpha^i h_i \left(v^i_x +\frac{1}{2}(w_x)^2\right)w_x\right] _x    + \gamma \beta h_3 (p_x w_x)_x=0
  \end{array}\right.
\end{eqnarray}
with  clamped boundary conditions at the left end $ \left.v^1, v^3, w, w_x, p~~\right|_{x=0}=0$ and the free boundary conditions at $x=L:$
\begin{eqnarray}
 \label{ivp}\left\{\hspace{-0.05in}
  \begin{array}{ll}
  %\left.v^1, v^3, w, w_x, p~~\right|_{x=0}=0,& \\
   \frac{1}{6} \alpha^2 h_2 v^3_x + +\alpha^1 h_1 (v^1_x+\frac{1}{2}(w_x)^2) + \frac{1}{3}\alpha^2 h_2v^1_{x}    +\frac{\rho_2h_2}{12}(2h_1-h_3) \ddot w-\frac{\alpha^2 h_2}{12} (2h_1 - h_3)w_{xx}=g^1(t),&\\
   \frac{1}{6} \alpha^2 h_2 v^1_x +\alpha^3 h_3 (v^3_x+\frac{1}{2}(w_x)^2) + \frac{1}{3}\alpha^2 h_2v^3_{x}   -\frac{\rho_2h_2}{12}(2h_3-h_1) \ddot w+\frac{\alpha^2 h_2}{12} (2h_3 - h_1)w_{xx}  -\gamma\beta h_3 p_x=g^3(t),&\\
   \beta h_3 p_x -\gamma \beta h_3 \left(v^3_x+\frac{1}{2}w_x^2\right) =-V(t),&\\
\frac{\alpha^1 (h_1)^3 + \alpha_3 (h_3)^3 + \alpha_2 h_2 ((h_1)^2 + (h_3)^2 -h_1 h_3)}{12} w_{xx}-\frac{\alpha_2 h_2}{12} (2h_1-h_3) v^1_{x}+\frac{\alpha_2 h_2}{12} (2h_3-h_1) v^3_{x}=M(t),\hspace{-0.2in}&\\
   &\\
   \frac{\rho_1 (h_1)^3 + \rho_3 (h_3)^3 + \rho_2 h_2 ((h_1)^2 + (h_3)^2 -h_1 h_3)}{12} \ddot w +G_2 H \phi^2  -\frac{\rho_2 h_2}{12} (2h_1-h_3) \ddot v^1+\frac{\rho_2 h_2}{12} (2h_3-h_1) \ddot v^1  &\\
   \quad +\frac{\alpha_2 h_2}{12} (2h_1-h_3) v^1_{xx} -\frac{\alpha_2 h_2}{12} (2h_3-h_1) v^3_{xx} +\sum_{i=1,3}\left[ \alpha^i h_i \left(v^i_x +\frac{1}{2}(w_x)^2\right)w_x\right] &\\
\quad   -\gamma\beta h_2 p_x w_x+ \frac{-\alpha^1 (h_1)^3 - \alpha_3 (h_3)^3 - \alpha_2 h_2 ((h_1)^2 + (h_3)^2 -h_1 h_3)}{12} w_{xxx}  =g(t),&\\
 %&\\
 (v^1, v^3, p, w, \dot v^1, \dot v^3, \dot p, \dot w)(x,0)=(v^1_0,  v^3_0, w_0, p_0,  v^1_1, v^3_1,  p_1, w_1).
  \end{array}\right.
\end{eqnarray}

%Note that damping terms  are intentionally discarded in  (\ref{perturbed}) to see the full effects of the controllers on the stabilization.

\subsection{Fully dynamic Rao-Nakra (R-N) beam model}
\label{modeling}
The model obtained above is still highly coupled. We assume that the viscoelastic layer is very thin and its stiffness is negligible, i.e. $\rho_2, \alpha^2\to 0$ as in  \cite{Hansen3}. This approximation retains the potential energy of shear and transverse kinetic energy so that the model above reduces to
\begin{eqnarray}
 \label{dbas} \left\{
  \begin{array}{ll}
 \rho_1h_1\ddot v^1-\alpha^1 h_1\left(v^1_{x} +\frac{1}{2}(w_x)^2\right)_x - G_2 \phi^2 = 0,   & \\
 \rho_3h_3\ddot v^3-\alpha^3 h_3\left(v^3_{x} +\frac{1}{2}(w_x)^2\right)_x+ G_2 \phi^2 + \gamma \beta h_3 p_{xx} =0,   & \\
 \mu  h_3 \ddot p   -\beta h_3   p_{xx}  + \gamma \beta h_3 \left(v^3_{x}+\frac{1}{2}(w_x)^2\right)_x= 0, &\\
  \rho \ddot w - K_1 \ddot{w}_{xx} + K_2 w_{xxxx} - G_2 H \phi^2_x-\gamma \beta h_3 p_{xx} -\left[\left(\sum\limits_{i=1,3}\alpha^i h_i \left(v^i_x +\frac{1}{2}(w_x)^2\right)\right)w_x\right]_x=0,&\\
 \phi^2=\frac{1}{h_2}\left(-v^1+v^3 + H w_x\right)&  \end{array} \right.
\end{eqnarray}
with initial conditions, and clamped boundary conditions at the left end $ [v^1,v^3,p,w,w_x](0)=0$ and the free boundary conditions at $x=L:$
\begin{eqnarray}
 \label{d-son} \left\{
  \begin{array}{ll}
% v^1(0)=v^3(0)= p(0)=0, &\\
 \alpha^1 h_1 (v^1_x+\frac{1}{2}(w_x)^2)=g^1(t), \quad \alpha^3 h_3 (v^3_x+\frac{1}{2}(w_x)^2) -\gamma \beta h_3 p_x=g^3(t),&\\
\beta h_3 p_x -\gamma \beta h_3 \left(v^3_x+\frac{1}{2}(w_x)^2\right) =-V(t),\quad K_2 w_{xx} = -M(t),&\\
  K_1 \ddot w_x -K_2 w_{xxx} + G_2 H \phi^2-\gamma\beta h_3 p_x(L) w_x(L)+\sum\limits_{i=1,3}\alpha^i h_i \left(v^i_x(L) +\frac{1}{2}(w_x)^2(L)\right)w_x(L)=g(t),&\\
(v^1, v^3, p, w, \dot v^1, \dot v^3, \dot p, \dot w)(x,0)=
(v^1_0,  v^3_0, p_0, w_0, v^1_1, v^3_1,  p_1, w_1)&
 \end{array} \right.
\end{eqnarray}
where $\rho=\rho_1h_1+ \rho_2h_2 + \rho_3 h_3,$ $K_1=\frac{\rho_1 h_1^3}{12} +\frac{\rho_3 h_3^3}{12},$ and $K_2=\frac{\alpha^1 h_1^3}{12}+\frac{\alpha^3 h_3^3}{12}.$

This model accounts for not only the coupling between the stretching and bending motions but also the dynamic behavior of the stretching of the piezoelectric layers, i.e. $\ddot v^1, \ddot v^3\ne 0.$ Note that the linearized model is first derived in \cite{Ozkan2}, and it is shown in \cite{Ozkan3} that certain vibration modes can not be controlled in the absence of $g^1(t).$

\subsection{Electrostatic (R-N) beam model}
\label{Nomag-RN}
By the electrostatic assumption, all magnetic effects are discarded, i.e. $\mu \to 0$ or $\mu\ddot p=0.$ Since the $p-$equation in \ref{dbas}) can be solved for $p,$  the system (\ref{dbas})-(\ref{d-son})  are simplified as the following
\begin{eqnarray}
 \label{d4-non} \left\{
  \begin{array}{ll}
   \rho_1h_1\ddot v^1-\alpha_{1}^1 h_1 \left(v^1_{x} +\frac{1}{2}(w_x)^2\right)_x - G_2  \phi^2 = 0,   & \\
    \rho_3h_3\ddot v^3-\alpha_{1}^3 h_3 \left(v^3_{x} +\frac{1}{2}(w_x)^2\right)_x+ G_2  \phi^2 = 0,   &\\
   \rho \ddot w - K_1 \ddot{w}_{xx} +K_2  w_{xxxx} - {G_2 H}  \phi^2_x-\left[\left(\sum\limits_{i=1,3}\alpha^i h_i \left(v^i_x +\frac{1}{2}(w_x)^2\right)\right)w_x\right]_x-\left(\gamma H_L w_x\right)_x V(t)=0, &
    %\phi^2=\frac{1}{h_2}\left(-v^1+v^3 + H w_x\right)&
  \end{array} \right.
\end{eqnarray}
with the initial conditions, clamped boundary condition at the left end $[v^1,v^3,w,w_x](0)=0,$ and free boundary conditions at $x=L:$
\begin{eqnarray}
\label{divp-non} \left\{
  \begin{array}{ll}
 %v^1(0)=v^3(0)= w(0)=w_x(0)=0, & \\
 \alpha_1^1 h_1 (v^1_{x} +\frac{1}{2}(w_x)^2)(L)= g^1(t), \quad  \alpha_1^3 h_3 (v^3_{x} +\frac{1}{2}(w_x)^2)(L)= -\gamma V(t),&\\
  K_2 w_{xx}(L) = -M(t),\quad  \left[ K_1 \ddot w_x-K_2 w_{xxx} + G_2 H \phi^2+ \left(\sum\limits_{i=1,3}\alpha^i h_i \left(v^i_x +\frac{1}{2}(w_x)^2\right)\right)w_x\right](L)=g(t),&\\
  (v^1, v^3, w, \dot v^1, \dot v^3, \dot w)(x,0)=  (v^1_0,v^3_0,w_0,v^1_1,v^3_1,w_1).
  \end{array} \right.
\end{eqnarray}
where we use (\ref{coef1}), and $g^3\equiv 0.$ Note that the voltage controller is the only controller for controlling the strains in the piezoelectric layer.

If we linearize (\ref{d4-non}) along the equilibrium $v^1=v^3=w=\dot v^1=\dot v^3=\dot w=0,$ or if we simply consider the linear stress-strain relationship in the beginning, we obtain the coupled system
\begin{eqnarray}
 \label{d4-lin}&& \left\{
  \begin{array}{ll}
   \rho_1h_1\ddot v^1-\alpha_{1}^1 h_1 v^1_{xx}- G_2  \phi^2 = 0,   & \\
    \rho_3h_3\ddot v^3-\alpha_{1}^3 h_3 v^3_{xx} + G_2  \phi^2 = 0,   &\\
   \rho \ddot w - K_1 \ddot{w}_{xx} +K_2  w_{xxxx} - {G_2 H}  \phi^2_x=0,
  \end{array} \right.\\
\label{divp-lin} &&\left\{
  \begin{array}{ll}
 v^1(0)=v^3(0)= w(0)=w_x(0)=0, \quad \alpha_1^1 h_1 (v^1_{x}= g^1(t), \quad \alpha_1^3 h_3 (v^3_{x} = -\gamma V(t),&\\
  K_2 w_{xx}(L) = -M(t),\quad  \left[ K_1 \ddot w_x-K_2 w_{xxx} + G_2 H \phi^2\right](L)=g(t),&\\
  (v^1, v^3, w, \dot v^1, \dot v^3, \dot w)(x,0)=  (v^1_0,v^3_0,w_0,v^1_1,v^3_1,w_1).
  \end{array} \right.
\end{eqnarray}
where the only coupling is due to the shear $\phi^2=\frac{1}{h_2}\left(-v^1+v^3 + H w_x\right)$ of the middle layer. This model is first obtained  in \cite{Ozkan2}. Later, it is shown to be exponentially stabilizable by the appropriate choice of controllers \cite{Ozkan3}. For elastic-viscoelastic-elastic model, there is currently a large literature for the controllability and stabilization  \cite{O-Hansen3,O-Hansen4} with various other boundary conditions.

\subsection{Fully dynamic  Mead-Marcus (M-M) beam model}
By the Mead-Marcus sandwich beam theory, we  assume that the longitudinal and rotational inertia terms, i.e. $\ddot v^1, \ddot v^3, \ddot w_{xx},$ are negligible in (\ref{perturbed})-(\ref{ivp}). Therefore, we obtain
\begin{eqnarray}
 \label{perturbeddd}\hspace{-0.05in} \left\{\hspace{-0.05in}
  \begin{array}{ll}
 - \frac{1}{6}\alpha^2 h_2 v^3_{xx}  -\left[\left(\alpha^1 h_1 (v^1_x+\frac{1}{2}(w_x)^2\right)+ \frac{1}{3}\alpha^2 h_2 v^1_{x}\right]_x  + \frac{\alpha^2 h_2}{12}(2h_1-h_3) w_{xxx}- G_2  \phi^2 = 0,&\\
- \frac{1}{6}\alpha^2 h_2 v^1_{xx}   -\left[\left(\alpha^3 h_3 (v^3_x+\frac{1}{2}(w_x)^2\right)+ \frac{1}{3}\alpha^2 h_2 v^3_{x}\right]_x    - \frac{\alpha^2 h_2}{12}(2h_3-h_1) w_{xxx} + G_2 \phi^2  + \gamma \beta h_3 p_{xx} = 0,   & \\
        \mu  h_3 \ddot p   -\beta h_3   p_{xx} + \gamma \beta h_3 \left(v^3_{x}+\frac{1}{2}(w_x)^2\right)_x= 0, &\\
  \rho \ddot w -  H G_2   \phi^2_x  + \frac{1}{12}\left(\alpha^1 h_1^3+\alpha^3 h_3^3 + \alpha^2 h_2(h_1^2+h_3^2 -h_1h_3) \right)w_{xxxx}  &\\
   \quad -\frac{\alpha^2 h_2}{12}(2h_1-h_3) v^1_{xxx} +\frac{\alpha^2 h_2}{12} (2h_3-h_1) v^3_{xxx} -  \left[\left(\sum\limits_{i=1,3}\alpha^i h_i\left(v^i_x +\frac{1}{2}w_x^2\right)- \gamma \beta h_3 p_x \right)w_x\right] _x=0,&\\
 %   &\\
%  \phi^2=\frac{1}{h_2}\left(-v^1+v^3 + H w_x\right),&
 \end{array} \right.
\end{eqnarray}
with clamped boundary conditions at the left end $  \left.v^1, v^3, w, w_x, p~~\right|_{x=0}=0$ and the free boundary conditions at $x=L:$
\begin{eqnarray}
 \label{ivpd}\hspace{-0.05in} \left\{\hspace{-0.05in}
  \begin{array}{ll}
 %\left.v^1, v^3, w, w_x, p~~\right|_{x=0}=0, \\
 %&\\
  \frac{1}{6}\alpha^2 h_2 v^3_{x}  +\left(\alpha^1 h_1 (v^1_x+\frac{1}{2}(w_x)^2\right)+ \frac{1}{3}\alpha^2 h_2 v^1_{x}   - \frac{\alpha^2 h_2}{12}(2h_1-h_3) w_{xx}=0,&\\
  \frac{1}{6}\alpha^2 h_2 v^1_{x}   + \alpha^3 h_3 \left(v^3_x+\frac{1}{2}(w_x)^2\right)+ \frac{1}{3}\alpha^2 h_2 v^3_{x}+ \frac{\alpha^2 h_2}{12}(2h_3-h_1) w_{xx}- \gamma \beta h_3 p_{x} = 0,\\
  \beta h_3 p_x -\gamma \beta h_3 \left(v^3_x+\frac{1}{2}(w_x)^2\right) =-V(t), \\
 \frac{\alpha^1 (h_1)^3 + \alpha_3 (h_3)^3 + \alpha_2 h_2 ((h_1)^2 + (h_3)^2 -h_1 h_3)}{12} w_{xx} -\frac{\alpha_2 h_2}{12} (2h_1-h_3) v^1_{x}+\frac{\alpha_2 h_2}{12} (2h_3-h_1) v^3_{x}=0,\\
     G_2 H \phi^2 +\frac{\alpha_2 h_2}{12} (2h_1-h_3) v^1_{xx} -\frac{\alpha_2 h_2}{12} (2h_3-h_1) v^3_{xx}   +\left(\sum\limits_{i=1,3}\left\{\alpha^i h_i\left(v^i_x +\frac{1}{2}(w_x)^2\right)\right\}+ \gamma \beta h_3 p_x \right)w_x+&\\
 \quad \frac{-\alpha^1 (h_1)^3 - \alpha_3 (h_3)^3 - \alpha_2 h_2 ((h_1)^2 + (h_3)^2 -h_1 h_3)}{12} w_{xxx} =g(t),&\\
 (v^1, v^3,  w, p, \dot v^1, \dot v^3, \dot w, \dot p)(x,0)=(v^1_0,  v^3_0, w_0,  p_0, v^1_1, v^3_1,   w_1, p_1).
 \end{array} \right. &&
\end{eqnarray}
Define the  constants by
\begin{eqnarray}
\nonumber \begin{array}{ll}
D:=\alpha^2 h_2^2 \left(4\alpha^1 h_1 + \alpha^2 h_2 + 4\alpha^3 h_3\right)+12 \alpha^1\alpha^3 h_1h_2 h_3,&\\
 A= \frac{1}{12}\left(\alpha^1h_1^3 + \alpha^3 h_3^3+ \frac{\alpha^2 h_2^2\left(3\alpha^2 h_2(\alpha^1h_1^3 + \alpha^3 h_3^3)+12\alpha^1\alpha^3h_1 h_2h_3 (h_1^2+h_3^2-h_1h_3)\right)}{D}\right),&\\
 B_1= \frac{\alpha^2\left[ \alpha^2 h_2^2 + 3 \alpha^1 h_1^2 +4\alpha^1 h_1 h_2 + 3 \alpha^3 h_3^2 + 4 \alpha^3 h_2 h_3\right]}{D} + 12 \frac{\alpha^1\alpha^3 H h_1 h_3}{D},~~  B_2= \frac{6\alpha^2 h_2 +12\alpha^1h_1 }{D},&\\
  %D_2= \frac{6\alpha^2 h_2 +12\alpha^3h_3 }{\alpha^2 h_2^2 \left(4\alpha^1 h_1 + \alpha^2 h_2 + 4\alpha^3 h_3\right)+12 \alpha^1\alpha^3 h_1 h_2h_3},&\\
 B_3=\frac{ \frac{1}{2}(\alpha^2)^2h_2^3h_3^2 + 2\alpha^1\alpha^2h_1h_2^2h_3^2-\alpha^1\alpha^2h_1^2h_2^2h_3}{D},~~B_4=\frac{(\alpha^2)^2 h_2^3 h_3  + 12 \alpha^1 \alpha^3_1 h_1h_2  h_3^2+ 4\alpha^1\alpha^2 h_1h_2^2 h_3 + 4\alpha^2\alpha^3_1 h_2^2h_3^2}{D},&\\
    B_5= \frac{6\alpha^1\alpha^2 h_1h_2^2 +(\alpha^2h_2)^2h_2 }{D},~~ C=\frac{12\left(\alpha^1h_1 +\alpha^2 h_2 +\alpha^3 h_3\right)}{D},~~ E_1=\frac{6\alpha^2 h_2 \left(\alpha^1 h_1 - \alpha^3 h_3 \right)}{D}, E_2= \frac{12\alpha^2 h_2h_3}{D},&\\
   E_3=-(\alpha^2 \alpha^3 h_2 h_3 +\alpha^1 \alpha^2 h_1 h_2+12\alpha^1 \alpha^3 h_1 h_3)C -\frac{(\alpha^1 h_1-\alpha^3 h_3)E_1}{D}
\end{array}
\end{eqnarray}
Consider $g(t)\equiv 0. $ Now we multiply the  first and second equations in (\ref{perturbeddd}) by $\left(\frac{1}{24}\alpha^2  h_2 + \frac{1}{12} \alpha^3 h_3\right) $ and \\$-\left(\frac{1}{24}\alpha^2  h_2 + \frac{1}{12} \alpha^1 h_1\right),$ respectively, and add these two equations. An alternate formulation is obtained as the following
\begin{eqnarray}
\nonumber&& \left\{ \begin{array}{ll}
  \rho \ddot w + A w_{xxxx} -  B_1G_2\gamma \beta  h_2h_3\phi^2_x + \gamma\beta B_3p_{xxx} +\frac{1}{C}(B_1E_1 w_x w_{xx} + E_1\phi^2_x w_x+ \gamma\beta  E_2 p_x w_x )_x&\\
  \quad + 3 E_3 (w_x)^2 w_{xx}=0,\\
  CG_2 \phi^2 -\phi^2_{xx} + B_1w_{xxx} +  B_2 p_{xx}+  E_1w_x w_{xx}=0,&\\
 \mu   \ddot p   + \beta \left(-B_4 p_{xx} + \gamma\beta h_2 h_3 G_2 B_2\phi^2 -\gamma  B_3 w_{xxx} +\gamma  B_5 w_x w_{xx}\right) =0,&
  \end{array} \right.\\
\label{perturbed-dumb}&& \left\{\begin{array}{ll}
  \left[ w, w_x, \phi^2, p~\right](0)=0,\quad  \left[A w_{xx} + \gamma\beta  B_3p_{x}+\frac{B_1E_1}{2C} (w_x)^2\right](L)=0,&\\
  \left[-B_4  p_{x}  -\gamma B_3 w_{xx}+ \frac{\gamma B_5}{2}(w_x)^2\right](L)= -\frac{V(t)}{h_3\beta},\quad \left[-\phi^2_{x} +B_1w_{xx} +  B_2 p_{x} +\frac{E_1}{2} (w_x)^2 \right](L)=0,&\\
   \left[A w_{xxx} - B_1G_2\gamma\beta  h_2 h_3 \phi^2 + \gamma\beta  B_3p_{xx}+\frac{1}{C}(B_1E_1 w_x w_{xx} + E_1\phi^2_x w_x+ \gamma\beta  E_2 p_x w_x )+ E_3 (w_x)^3\right](L)=0&\\
( w, p, \dot w, \dot p)(x,0)=(w_0,  p_0, w^1, p^1).
\end{array}\right.
\end{eqnarray}
%In fact, the boundary conditions can be further simplified to obtain
%\begin{eqnarray}
%\nonumber && \left[ w, w_x, \phi^2, p\right](0)=0, \quad \left[w_{xx}, \phi^2_x, p_x\right](L)=0, \\
%\nonumber &&\left[- A w_{xxx} + B_1G_2 \gamma \beta h_2 h_3 \phi^2 - \gamma\beta B_3p_{xx}\right](L)=0.
%\end{eqnarray}

\subsection{Electrostatic (M-M) beam model}
%By using (\ref{coef1}), define
%\begin{eqnarray}
%\nonumber \begin{array}{ll}
 %\tilde A= \frac{1}{12}\left(\alpha^1h_1^3 + \alpha^3 h_3^3\right.&\\
 % \left.+ \frac{\alpha^2 h_2^2\left(3\alpha^2 h_2(\alpha^1h_1^3 + \alpha^3_1 h_3^3)+12\alpha^1\alpha^3_1h_1 h_2h_3 (h_1^2+h_3^2-h_1h_3)\right)}{D}\right),&\\
 %\tilde D=D(\alpha^3 \rightarrow \alpha^3_1),~\tilde C=C(\alpha^3 \rightarrow \alpha^3_1)&\\
 %\tilde B_i= B_i(\alpha^3 \rightarrow \alpha^3_1), \quad i=1,\ldots,5&\\
  % \tilde E_i=E_i(\alpha^3 \rightarrow \alpha^3_1), \quad i=1,\ldots,2.&
%\end{array}
%\end{eqnarray}
Notice that if the magnetic effects in (\ref{perturbed-dumb}) are neglected, i.e. $\ddot \mu p\equiv 0,$ the last equation can be solved for $p_{xx}.$ Then we obtain the following model with simplified boundary conditions
 \begin{eqnarray}
\label{abstractMM}&&\left\{ \begin{array}{ll}
 \rho \ddot w + {\tilde A} w_{xxxx} - G_2\gamma \beta  h_2h_3{\tilde B}\phi^2_x   +(N_1 w_x w_{xx} + N_2w_x\phi^2_x +  N_3 (w_x)^3 +\frac{\gamma E_2 }{h_3 \tilde C B_4} V(t)H_Lw_x )_x=0, &\\
G_2{ {\tilde C}} \phi^2 -\phi^2_{xx} + {{\tilde  B}} w_{xxx}+N_4 w_x w_{xx}=0,
\end{array}\right.\\
\label{MM-bc} &&\left\{\begin{array}{ll}
  \left[ w, w_x, \phi^2~\right](0)=0,\quad  \left[-\phi^2_{x} +\tilde B w_{xx} +  \frac{1}{2}N_4 (w_x)^2 \right](L)=-\frac{B_2V(t)}{\beta h_3 B_4},&\\
   \left[\tilde A w_{xx} +\frac{1}{2}\left(N_1+ \frac{\gamma^2\beta B_3E_2}{\tilde C B_4}\right) (w_x)^2\right](L)=-\frac{\gamma B_3V(t)}{h_3 B_4},&\\
   \left[\tilde A w_{xxx} -  G_2\gamma \beta  h_2h_3{\tilde B} \phi^2 +N_1 w_x w_{xx} + N_2w_x\phi^2_x +  N_3 (w_x)^3 +\frac{\gamma E_2 }{h_3 \tilde C B_4} V(t)w_x\right](L)=0,&\\
( w, \dot w)(x,0)=(w_0,  w^1)
\end{array}\right.
\end{eqnarray}
where the coefficients are functions of material parameters
\begin{eqnarray}
\nonumber \begin{array}{ll}
\tilde  A= A-\frac{\gamma^2\beta B_3^2}{B_4},~~ \tilde B=B_1 -\frac{\gamma\beta B_2 B_3}{B_4},~~ \tilde B_3=B_3 + \frac{\gamma\beta h_2 h_3 B_1 B_2}{C},~~ \tilde B_4=B_4 + \frac{\gamma\beta h_2 h_3 B_2^2}{C},~~ \tilde B_5=B_5 - \frac{\beta h_2 h_3 B_2 E_1}{C}&\\
\tilde C=\frac{C\tilde B_4}{B_4},~~  N_1=\frac{B_1 E_1 }{C}+\gamma^2\beta B_3\left(\frac{  B_5}{B_4}- \frac{E_2}{\tilde C B_4}\right),~~ N_2=\frac{E_1}{C}+ \frac{\gamma^2\beta^2 h_2 h_3 B_2 E_2}{C\tilde C B_4},~~ N_3=E_3 + \frac{\gamma^2\beta E_2 \tilde B_5}{6\tilde C B_4},&\\
N_4=E_1+ \frac{\gamma B_2 B_5}{B_4}.&
\end{array}
\end{eqnarray}
 Unlike the electrostatic Rao-Nakta model (\ref{d4-non}), the voltage control simultaneously controls the shear of the middle layer and the bending motion of the composite.
%are ${\tilde A}= A-\frac{\gamma^2 \beta h_2 h_3\tilde B_3^2}{\tilde B_4}>0, {\tilde B}=B_1-\frac{\gamma^2 h_3\beta\tilde  B_2 \tilde B_3}{\tilde B_4}>0, {{\tilde C}}= C+ \frac{\gamma^2 \beta h_3 \tilde B_2^2}{\tilde B_4},$ ${\tilde E_1}= E_1+\frac{\gamma^2 \beta h_3 \tilde B_2 B_5}{\tilde B_4},$ $N_1,
%\ldots N_4$ are functions of all other constants.
This model is similar to the model obtained in \cite{B-R} except that the inclusion of longitudinal strains is more rigorous here. Note that If we linearize (\ref{abstractMM}) along the equilibrium $\phi^2=w=\dot w=0,$ or if we simply consider the linear stress-strain relationship in the beginning, we obtain the coupled system
 \begin{eqnarray}
\label{abstractMM-lin}&&\left\{ \begin{array}{ll}
 \rho \ddot w + {\tilde A} w_{xxxx} - G_2\gamma \beta  h_2h_3{\tilde B}\phi^2_x=0, &\\
G_2{ {\tilde C}} \phi^2 -\phi^2_{xx} + {{\tilde  B}} w_{xxx}=0,&\\
\end{array}\right.\\
\label{MM-bc-lin} &&\left\{\begin{array}{ll}
  \left[ w, w_x, \phi^2, p~\right](0)=0,\quad    \left[-\phi^2_{x} +\tilde B w_{xx}\right](L) =-\frac{B_2V(t)}{\beta h_3 B_4},&\\
  \tilde A w_{xx} (L)=-\frac{\gamma B_3V(t)}{h_3 B_4},\quad  \left[\tilde A w_{xxx} -  G_2\gamma \beta  h_2h_3{\tilde B} \phi^2 \right](L)=0,&\\
( w, \dot w)(x,0)=(w_0, w^1).
\end{array}\right.
\end{eqnarray}
The model (\ref{abstractMM-lin})-(\ref{MM-bc-lin}) has been recently proved to be exponentially stabilizable by only the voltage controller \cite{Ozkan6}.

\begin{rmk}
The fully dynamic composite beam models (\ref{dbas}) and (\ref{perturbed-dumb})  with appropriate choice of states {\bf y} and input $u(t)$ can be formulated in the following form
$\dot {\bf y}=(\mc A +\mc N){\bf y} + \mc B u(t)$ where $\mc B$ is a boundary control operator  for which the observation $\mc B^*{\bf y}$ corresponds to the total induced current accumulated at the electrodes for the voltage input  \cite{Ozkan3}. This is completely electro-magnetic.

Unlike fully dynamic models, electrostatic models  (\ref{d4-non}) and (\ref{abstractMM}) with appropriate choice of states {\bf y} and input $u(t)$ can be formulated  $\dot {\bf y}=(\mc A +\mc N){\bf y} + (\mc B_1 + \mc B_2{\bf y}) u(t)$ where $\mc B_1$ is a boundary control operator and $\mc B_2$ is an unbounded bilinear control operator. This form is similar to the one obtained in (\ref{s-s}).
\end{rmk}
\begin{rmk}
For all models (\ref{dbas}),(\ref{d4-non}), (\ref{perturbed-dumb}) and (\ref{abstractMM}), a shear type of damping (due to the viscoelastic layer) can be added by replacing the term $G_2 \phi^2$ by $G_2 \phi^2 + \tilde G_2 \dot \phi^2$ where $\tilde G_2$ is the damping coefficient  \cite{Hansen3}.
\end{rmk}
\begin{rmk} \label{har} For energy harvesting applications, the nonlinear models (\ref{d4-non},\ref{abstractMM}) and the linear models (\ref{d4-lin},\ref{abstractMM-lin}) are coupled to the extra circuit equation
\begin{eqnarray}
\frac{2\alpha_{11}\ep3 b L}{\alpha_1h} \dot V(t)+ i(t) =-\gamma_3 (h_1+h_2) b \dot w_x(L,t) 
\end{eqnarray}
where $i(t)=\frac{1}{R}V(t)$ is the current generated by the piezoelectric beam, $R$ is the resistance of the attached circuit, see \cite{E-Inman,Shubov} for two-layer counterparts.
\end{rmk}
\section{Preliminary stabilization results for the semi-discrete approximations of a single piezoelectric beam}
\label{Num}
For the numerical analysis of the piezoelectric beam models, the Galerkin-based Finite Element Method (FEM) is a very standard way to see the effectiveness of the control design on the first several vibrational modes of the composite. However, our models are highly coupled and require more careful treatment for the spurious high frequency modes which may cause the so-called spill-over effect.  There are recent attempts to mathematically address these issues  such as, filtering techniques  for the high frequencies in finite differences \cite{Wehbe,Alabou,M-Z}, multigrid techniques \cite{N-Z} or mixed finite element methods \cite{A-C,Ito,M-Z}. Among these, the semi-discretized Finite Difference technique with the introduction of an artificial (numerical) viscosity term is performing well for many de-coupled equations, see \cite{T-Z,Leon}.The main idea to suppress the spurious high frequency oscillations is either (i) to control the projected high-frequency solutions on corresponding subspace so these spurious solutions are filtered, or (ii)  we add an extra boundary control. In fact, in a recent paper \cite{Ozkan6}, the filtered Finite Difference Method has powerful results for the linearized (M-M) beam equations (\ref{perturbed-dumb}) with PID controllers.

 The aim of this section is to present some simple numerical experiments in order to show that  $\mc B^*-$type stabilizing feedback controller can be designed in a similar fashion as for the linearized models.  The fully dynamic models are currently left out to save space. We only provide preliminary stabilization results for the electrostatic (E-B) and (M-T) single beam cases. Furthermore, we skip the convergency and consistency analyses. These will be comprehensively discussed in \cite{K-O}. The total energy of a single beam for (\ref{homo-NoM-EB}) and (\ref{homo-NoM-MT}) is
\begin{eqnarray}
  \nonumber &&\mb{E}(t)=\frac{h}{2} \left\{
  \begin{array}{ll}
    \int_0^L \left[\rho \left(\dot v^2+\frac{ h^2}{12} \dot w_x^2 + \dot w^2\right) +    \left( \alpha_{11} \left(v_x + \frac{1}{2}w_x^2\right)^2+ \frac{\alpha_1 h^2}{12}w_{xx}^2 \right)\right]~dx, & \text{(E-B)}\\
\int_0^L \left[\rho \left(\dot v^2+\frac{ h^2}{12} \dot \psi^2 + \dot w^2\right)+  \left(\alpha_{11} \left(v_x + \frac{1}{2}w_x^2\right)^2+ \alpha_1\frac{ h^2}{12}\psi_x^2 \right)+\alpha_3(w_x+\psi)^2\right]~dx,& \text{(M-T)}\\
  \end{array}
\right.
\end{eqnarray}
For $c_i>0$, we design a $\mc B^*-$type stabilizing control law as in Table \ref{fcontrols}.
 \bgroup
\def\arraystretch{1.2}
\begin{table}[h]
\centering
\begin{tabular}{|c|c|} \hline
Euler-Bernoulli (E-B) & Mindlin-Timoshenko (M-T) \\ \hline\hline
$ V(t)=c_1 \left(\dot v(L,t)+ \int_0^L w_x \dot w_{x} ~dx\right) $ & $ V(t)=c_4 \left(\dot v(L,t)+ \int_0^L w_x \dot w_{x} ~dx\right) $  \\ \hline
$m(t)=-c_2 \dot w_x(L,t)$ & $m(t)=-c_5 \dot \psi(L,t)$   \\ \hline
$g(t)=c_3 \dot w_t(L,t)$ & $g(t)=-c_6 \dot w_t(L,t)$  \\ \hline
\end{tabular}
\caption{\small Stabilizing feedback controllers. Notice that the voltage controller $V(t)$ has the nonlinear term $\int_0^L w_x \dot w_x dx.$ This is the contribution of nonlinearity  to the $\mc B^*-$feedback law.}
\label{fcontrols}
\end{table}
Notice that, $\mc B^*-$feedback  has a nonlinear integral controller $\int_0^L w_x \dot w_x dx,$ which is the contribution of the nonlinearity. The main advantage of an integral controller is that it provides high-gain feedback at low frequencies, and therefore, integral controllers can overcome creep
and hysteresis effects and lead to precision positioning (since the vibrational dynamics is not dominant at low frequencies) \cite{Review}. With this choice,  the energy for each model is dissipative
\begin{eqnarray}
  \nonumber &&\frac{d\mb{E}(t)}{dt}=\left\{
  \begin{array}{ll}
  -c_1  |\dot v(L,t)+ \int_0^L w_x \dot w_{x} ~dx| ^2-c_2 |\dot w_x(L,t)|^2-c_3 |\dot w(L,t)|^2   & \text{(E-B)}\\
-c_1  |\dot v(L,t)+ \int_0^L w_x \dot w_{x} ~dx| ^2-c_2 |\dot \psi(L,t)|^2-c_3 |\dot w(L,t)|^2 &  \text{(M-T)}\\
  \end{array}\le 0.
\right.
\end{eqnarray}
Now consider the  discretization of the interval $[0,L]$ with the fictitious points $x_{-1}$ and $x_{N+1}$ as the following
$$x_{-1}<0=x_0< x_1 <x_2, \ldots < x_N = L < x_{N+1}, \quad x_i=i\cdot dx, \quad i=-1,0,1,\ldots, N, N+1, ~~dx=\frac{L}{N+1}.$$ The following are the second order finite difference approximations for different order derivatives:
\begin{eqnarray}
\begin{array}{ll}
z_{x}=\frac{z(x_{i+1})-z(x_{i-1})}{2dx}+O(dx^2), \quad {\rm or}\quad \frac{3z(x_i)-4z({x_{i-1}})+z({x_{i-2}})}{2dx}+O(dx^2),&\\
\quad z_{xx}=\frac{z(x_{i+1})-2z(x_i)+z(x_{i-1})}{dx^2}+O(dx^2),\quad z_{xxxx}=\frac{z(x_{i+2})-2z(x_{i+1})+6z(x_{i})-4z (x_{i-1})+z(x_{i-2})}{dx^4}+O(dx^2).
\end{array}
\end{eqnarray}
Henceforth, to simplify the notation, we use $z(x_{i})=z_i.$
We also non-dimensionalize the time variable  $t=A_1 t^*~ {\rm with} ~t^*\to t.$ %where $A_1=L\sqrt{\frac{\rho}{\alpha_{11}}}\sim 0.023.$

We consider a  sample piezoelectric beam with height  $L=1 {\rm  m},$ and thickness $h= 0.01{\rm m}.$ The  material constants are $\rho=7600$ kg/{m$^3$},   $\alpha_1=1.4\times 10^7$ N/{m{$^2$}}, $\alpha_3=4.5\times 10^5$ N/{m{$^2$}}, $\gamma=10^{-3}$ C/{m$^2$}, $\beta= 10^{6} {\rm m/F},$   $G_2=100$ GN/${\rm m}^2$. We consider $N=60$ with the initial data $w(x,0)= v(x,0)=\dot v(x,0)= 10^{-3} e^{\left(\frac{x-0.5L}{0.1L}\right)^2}$ and $\dot w(x,0)=\psi(x,0)=\dot\psi(x,0)=0.$ The simulations are computed for (non-dimensionalized) the  final time $T_{\rm final}=300$.

The filtering technique is successfully applied in \cite{T-Z,Roventa} so that high frequency solutions, causing artificial instability in the approximated solution,  are filtered by adding a viscosity term $(dx)^2 u_{xxt}$ to the wave equation and $ w_{xxt}$ to the fourth order  beam equation. Both terms vanish uniformly as $dx\to 0.$ We adopt the same idea to derive a filtered semi-discrete scheme in finite differences to simulate the effects of the stabilizing controller.  For this purpose, additional viscosity terms (boxed terms in the following) $(dx)^2 v_{xxt}, \frac{1}{2} w_{xxt}, (dx)^2\psi_{xxt}$ to the (E-B) equations and  $(dx)^2 v_{xxt}, (dx)^2 w_{xxt}, (dx)^2\psi_{xxt}$ to the (M-T) equations  are added with appropriate coefficients:
\begin{eqnarray}
  \label{dis-es-EB}
&& (E-B) \left\{
\begin{array}{ll}
   \ddot v_i-\boxed{\frac{\dot v_{i+1}-2\dot v_i+\dot v_{i-1}}{dx^2}}- \left( \frac{v_{i+1}-2v_i+v_{i-1}}{dx^2} + \frac{1}{L}\frac{w_{i+1}-w_{i-1}}{2dx}  \frac{w_{i+1}-2w_i+w_{i-1}}{dx^2}\right) = 0   & \\
   \ddot w -\frac{h^2}{12L^2}\frac{\ddot w_{i+1}-2\ddot w_i+\ddot w_{i-1}}{dx^2}-\frac{1}{2}\frac{h^2}{12}\boxed{\frac{\dot w_{i+1}-2\dot w_i+\dot w_{i-1}}{dx^2}}+\frac{h^2}{12}\frac{w_{i+2}-4w_{i+1}+6w_i-4w_{i-1}+w_{i-2}}{dx^4}&\\
~ - \frac{1}{L}\left[ \frac{v_{i+1}-v_{i-1}}{2dx} \frac{w_{i+1}-2w_i+w_{i-1}}{dx^2} +  \frac{v_{i+1}-2v_i+v_{i-1}}{dx^2}  \frac{w_{i+1}-w_{i-1}}{2dx}  + \frac{3}{2L^3} \left[\frac{w_{i+1}-w_{i-1}}{2dx}\right]^2  \frac{w_{i+1}-2w_i+w_{i-1}}{dx^2} \right]&\\
 ~ -\frac{ \gamma}{\alpha_{11} h} \left(\frac{w_{i+1}-2w_i+w_{i-1}}{dx^2}\right)V(t)=0, \quad i=1,\ldots, N-1&\\
 &\\
 v_0=w_0=0, w_{-1}=w_1, \quad  ~\frac{h^2}{12L^2}\frac{w_{N+1}-2w_N+w_{N-1}}{dx^2}=-m(t)&\\
  \frac{3v_N-4v_{N-1}+v_{N-2}}{2dx} +  \frac{1}{2}\left(\frac{3w_N-4w_{N-1}+w_{N-2}}{2dx}\right)^2 = -\frac{\gamma}{\alpha_{11} h}V(t)&\\
\frac{3\ddot w_N-4\ddot w_{N-1}-\ddot w_{N-2}}{2dx}+\frac{h^2}{12L^2}\frac{2w_{N+1}-5w_N+2w_{N-1} + 4w_{N-2}-4w_{N-3}}{2dx^3} =g(t),&
\end{array}\right.
\end{eqnarray}
\begin{eqnarray}
 \label{dis-es-MT}
&&(M-T)   \left\{
\begin{array}{ll}
   \ddot v_i-\boxed{\frac{\dot v_{i+1}-2\dot v_i+\dot v_{i-1}}{dx^2}}- \left( \frac{v_{i+1}-2v_i+v_{i-1}}{dx^2} + \frac{1}{L}\frac{w_{i+1}-w_{i-1}}{2dx}  \frac{w_{i+1}-2w_i+w_{i-1}}{dx^2}\right) = 0,   & \\
   \ddot \psi_i -\frac{\alpha_1}{\alpha_{11}}\boxed{\frac{\dot \psi_{i+1}-2\dot \psi_i+\dot \psi_{i-1}}{dx^2}}-\frac{\alpha_1}{\alpha_{11}}\frac{ \psi_{i+1}-2 \psi_i+\psi_{i-1}}{dx^2}+\frac{12\alpha_3}{ \alpha_{11} h^2} \left(L\frac{w_{i+1}-w_{i-1}}{2dx} +L^2 \psi_i\right)=0, &\\
   \ddot w -\frac{\alpha_3}{\alpha_{11}}\boxed{\frac{\dot w_{i+1}-2\dot w_i+\dot w_{i-1}}{dx^2}} -\frac{\alpha_3}{\alpha_{11}} \left( \frac{w_{i+1}-2w_i+w_{i-1}}{dx^2}  +\frac{\psi_{i+1}-\psi_{i-1}}{2dx}\right)&\\
      ~  - \frac{1}{L}\left[ \frac{v_{i+1}-v_{i-1}}{2dx} \frac{w_{i+1}-2w_i+w_{i-1}}{dx^2} +  \frac{v_{i+1}-2v_i+v_{i-1}}{dx^2}  \frac{w_{i+1}-w_{i-1}}{2dx}  + \frac{3}{2L^3} \left[\frac{w_{i+1}-w_{i-1}}{2dx}\right]^2  \frac{w_{i+1}-2w_i+w_{i-1}}{dx^2} \right]&\\
  ~ - \frac{ \gamma}{\alpha_{11} h}\left(\frac{w_{i+1}-2w_i+w_{i-1}}{dx^2}\right)V(t)=0, \quad i=1,\ldots, N-1&\\
&\\
  v_0=\psi_0=w_0=0, \quad  \frac{3v_N-4v_{N-1}+v_{N-2}}{2dx} +  \frac{1}{2}\left(\frac{3w_N-4w_{N-1}+w_{N-2}}{2dx}\right)^2  = -\frac{\gamma}{\alpha_{11} h}V(t)&\\
    %\left[\beta_{3} h p'-{\gamma_3\beta_3}h  \left(v' + \frac{1}{2}w'^2\right)=-V(t)\right]_{x=L} &\\
      \alpha_3 h \frac{3\psi_N-4\psi_{N-1}+\psi_{N-2}}{2dx}  =m(t),\quad  \alpha_3 h \left(\frac{3w_N-4w_{N-1}+w_{N-2}}{2dx}+\psi_N\right)=g(t)
\end{array}\right.
\end{eqnarray}
where the approximated controllers are designed
\begin{eqnarray}
\nonumber &&{\rm (E-B)} \left\{
  \begin{array}{ll}
    V(t) =c_1 \left(\dot v_{N} + \frac{dx}{3}\frac{d}{dt}\left(\left(\frac{3w_N-4w_{N-1}+w_{N-2}}{2dx}\right)^2+4\sum\limits_{i=1}^{N/2} \left(\frac{w_{2i}-w_{2i-2}}{2dx}\right)^2 \right.\right.&\\
    \left.\left.\quad\quad+ 2\sum\limits_{i=1}^{N/2-1} \left(\frac{w_{2i+1}-w_{2i-1}}{2dx}\right)^2\right) \right), &\\
m(t)=-c_2\frac{d}{dt} \left(\frac{3w_N-4w_{N-1}+w_{N-2}}{2dx}\right), \quad g(t)=c_3 \dot w_N &
  \end{array}
\right.\\
\label{controllers} &&{\rm (M-T)} \left\{
  \begin{array}{ll}
    V(t) =c_4 \left(\dot v_{N} + \frac{dx}{3}\frac{d}{dt}\left(\left(\frac{3w_N-4w_{N-1}+w_{N-2}}{2dx}\right)^2+4\sum\limits_{i=1}^{N/2} \left(\frac{w_{2i}-w_{2i-2}}{2dx}\right)^2 \right.\right.&\\
    \left.\left.\quad\quad+ 2\sum\limits_{i=1}^{N/2-1} \left(\frac{w_{2i+1}-w_{2i-1}}{2dx}\right)^2\right) \right), &\\
m(t)=-c_5\dot \psi_N, \quad g(t)=-c_6 \dot w_N &
  \end{array}
\right.
\end{eqnarray}
where the feedback gains $c_1,\ldots c_6>0$ are chosen appropriately in the numerical code. The Differential Algebraic system of equations are solved by the Mathematica's {\bf NDSolve} command with AccuracyGoal and PrecisionGoal options to be set to MachinePrecision/2. The simulations are run for three major cases:
\begin{itemize} \item  Fully controlled : $V(t), m(t),g(t)\ne 0,$
\item  Partially controlled: $V(t)\equiv 0 $ and $m(t),g(t)\ne 0,$
\item No control: $V(t), m(t),g(t)\equiv  0.$
\end{itemize}

%The voltage controller $V(t)$ in (\ref{controllers})  stabilizes the stretching solutions faster in comparison to the controllers $m(t), g(t)$ stabilizing the bending or shear solutions, see Tables \ref{g1} and \ref{g2}. We also considered the partially controlled case where. $V(t)\equiv 0, m(t), g(t)\ne 0$. With only two controllers, stretching motions are not controlled at all. However,  the polynomial decay for the bending and shear solutions are still recovered.

In the fully controlled case, the feedback controller corresponding to the voltage control $V(t)$ in (\ref{controllers}) is very powerful to  exponentially stabilize the stretching solutions $v(x,t)$ for both models, see Tables \ref{g1} and \ref{g2}.   In the partially controlled case for both  (E-B) and (M-T) models, two controllers $m(t)$ and $g(t)$ are able to decay the  bending and shear to zero polynomially. Since the voltage controller is freed in this case, the stretching solutions in both (E-B) and (M-T) cases are not controlled at all. All results are also compared to the uncontrolled models in Tables \ref{g1} and \ref{g2} .

The tip velocities $\dot v_N,  \dot w_N$ for the (E-B) model and   $\dot v_N, \dot \psi_N, \dot w_N$ for the (M-T) model  are simulated in Tables \ref{g4}  and \ref{g5}, respectively. It can be observed that the $\dot v_N$ decay to zero exponentially for both models. The polynomial decay pattern for the $\dot \psi_N$ and $\dot w_N$ is still observed.
The energy of the fully controlled models, shown  in Table \ref{g3},  decay to zero in a few seconds. Unlike the fully elastic models \cite{Lag-L}, it is our observation that the voltage controller $V(t)$, controlling stretching equations, dominates the other controls in the approximated schemes (\ref{dis-es-EB}) and (\ref{dis-es-MT}).
This implies that the boundary control design for the approximated solutions (\ref{dis-es-EB}) and (\ref{dis-es-MT}) should  be improved by maybe adding additional viscosity terms (boundary type or distributed type). Adaptive controllers may also be considered to improve the decay rate of bending and shear solutions. Analytical work  to find out the decay rates of the nonlinear models is currently under investigation \cite{Ozkan7}.
 \bgroup
\def\arraystretch{1.0}
\begin{table}[htb]
\begin{tabular}{cc}
\centering Stretching - $v(x,t)$ & Bending - $w(x,t)$ \\
\includegraphics[width=3.2in]{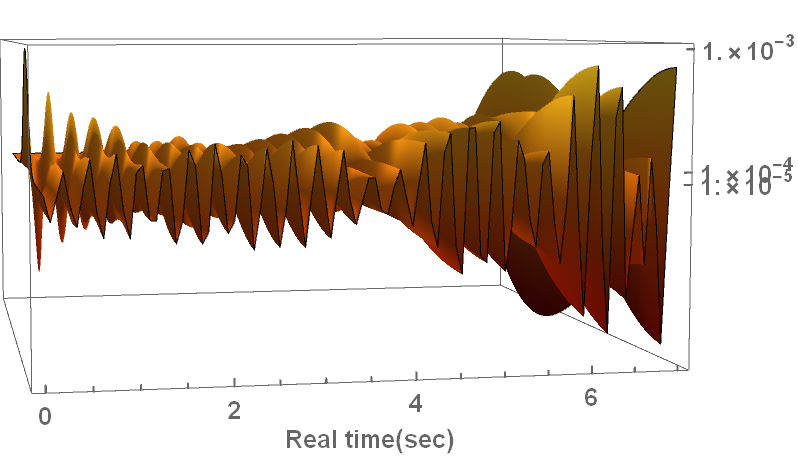} & \includegraphics[width=3.2in]{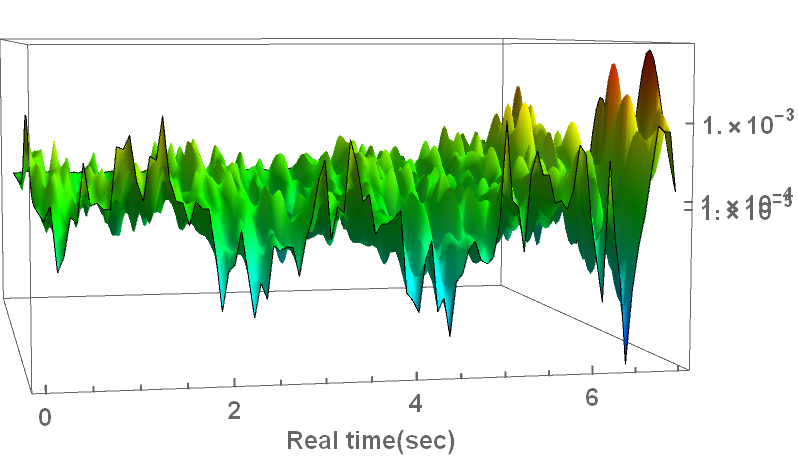}  \\
\includegraphics[width=3.2in]{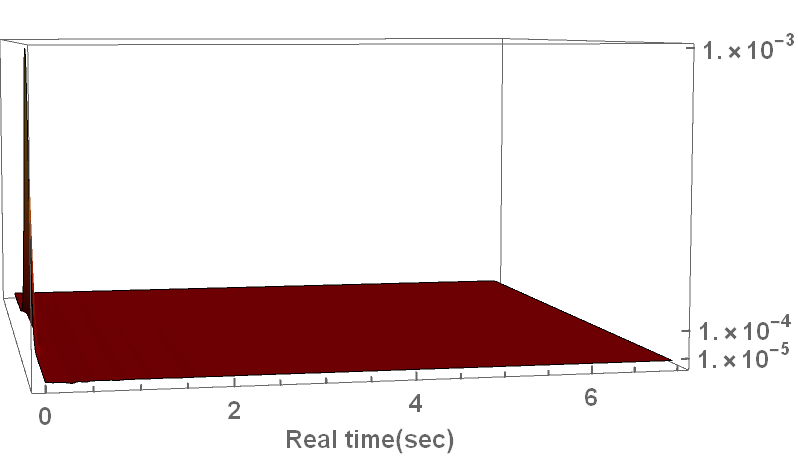} & \includegraphics[width=3.2in]{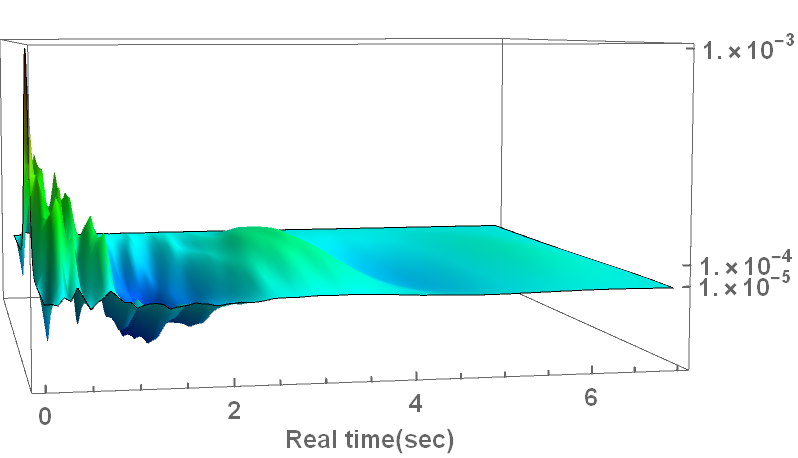}  \\
\includegraphics[width=3.2in]{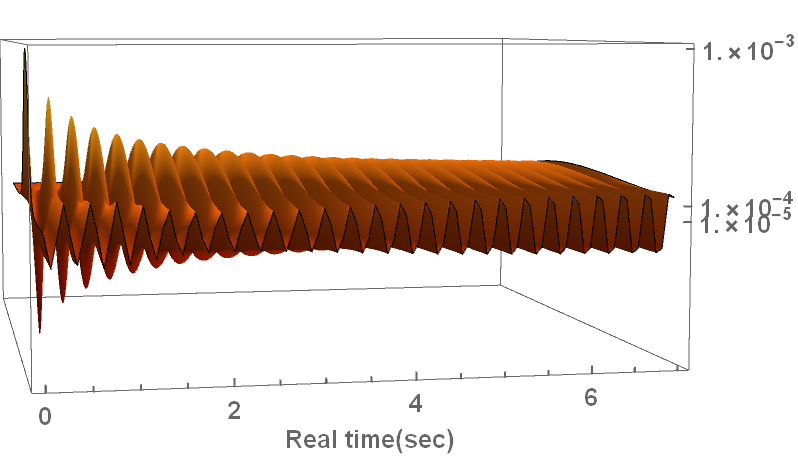} & \includegraphics[width=3.2in]{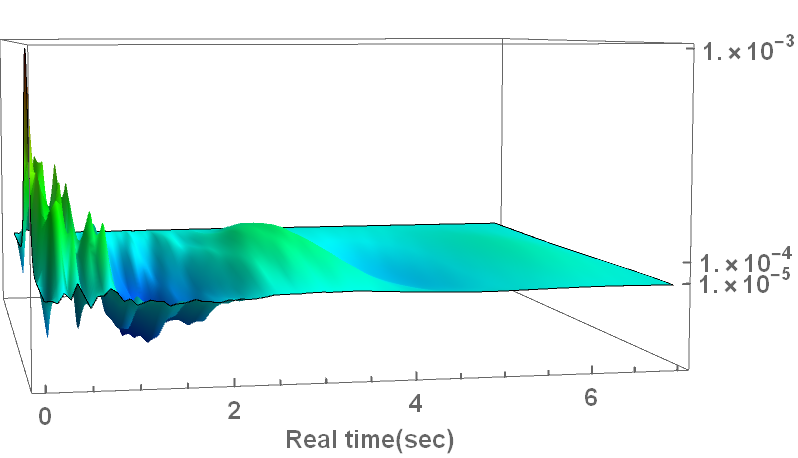}  \\ \hline
\end{tabular}
\caption{\footnotesize Graphs in each column show simulations of the approximated solutions $\{v, w\}$ of the (E-B) closed-loop system for uncontrolled, controlled, and partially controlled ($V(t)=0$) cases, respectively. The $t-$axis is real-time. }
\label{g1}
\end{table}

 \bgroup
\def\arraystretch{1.1}
\begin{table}[h]
\begin{tabular}{ccc}
\centering Stretching - $v(x,t)$ & Bending - $w(x,t)$ & Rotation angle $\psi_N(x,t)$\\
\includegraphics[width=2.1in]{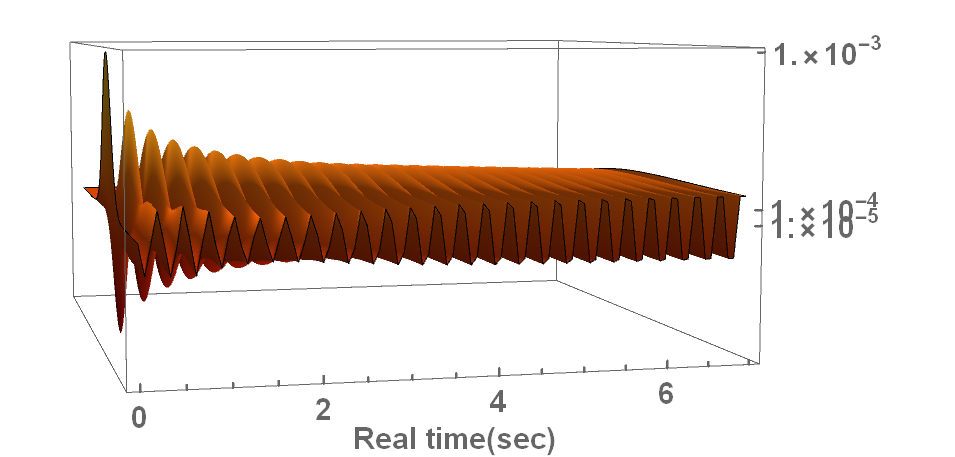} &\includegraphics[width=2.1in]{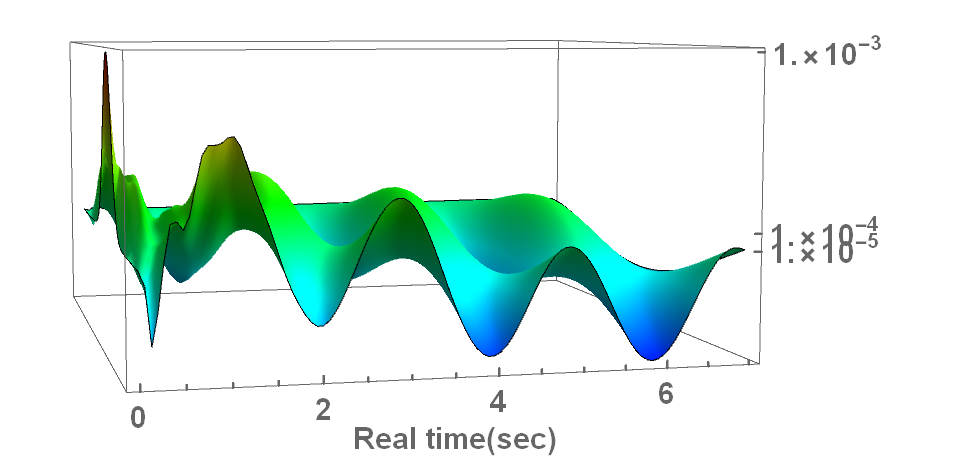} &\includegraphics[width=2.1in]{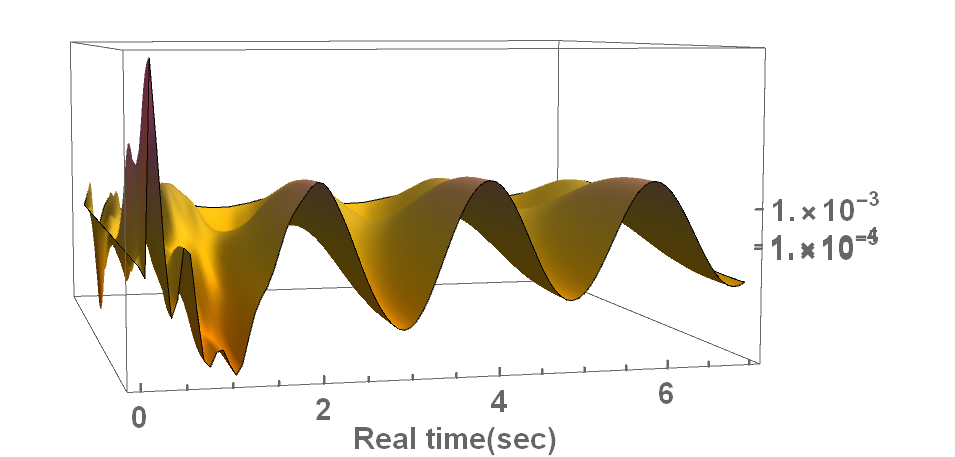}  \\
\includegraphics[width=2.1in]{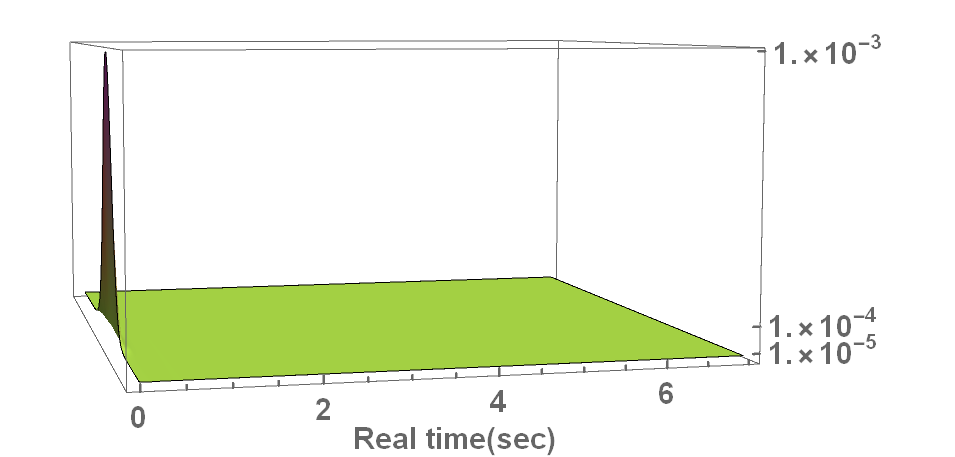} &\includegraphics[width=2.1in]{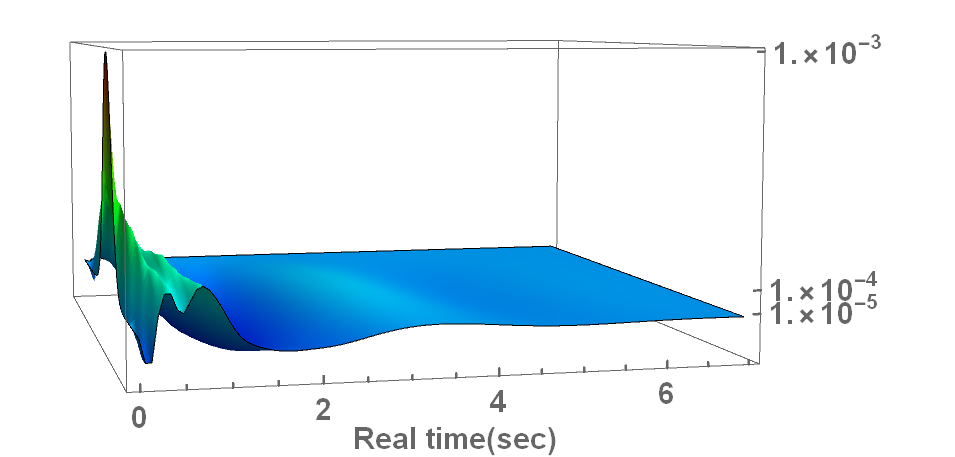}  & \includegraphics[width=2.1in]{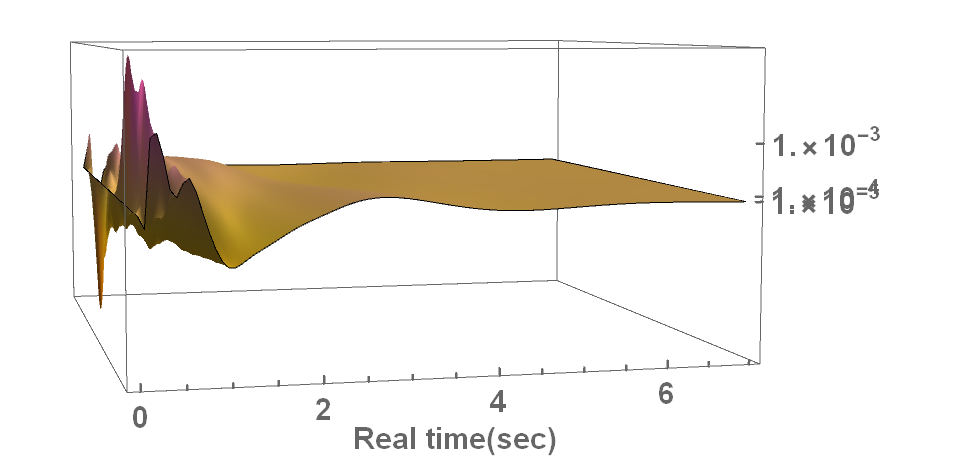}   \\
\includegraphics[width=2.1in]{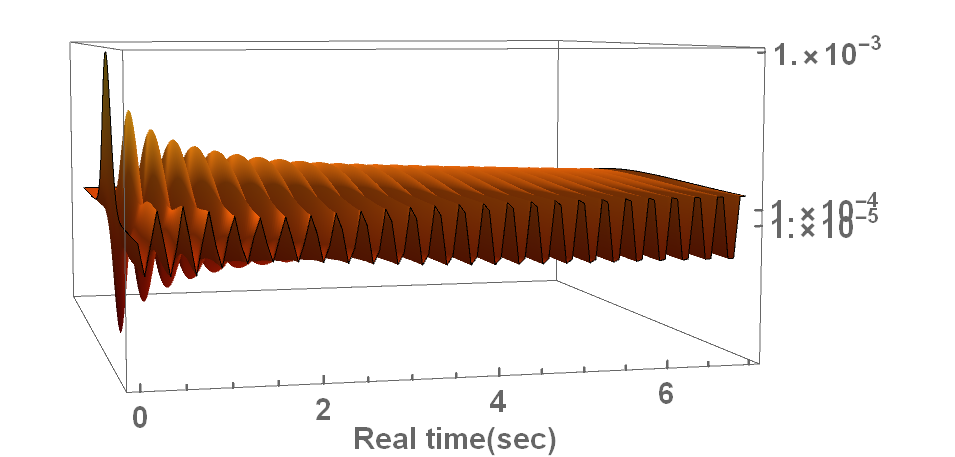} & \includegraphics[width=2.1in]{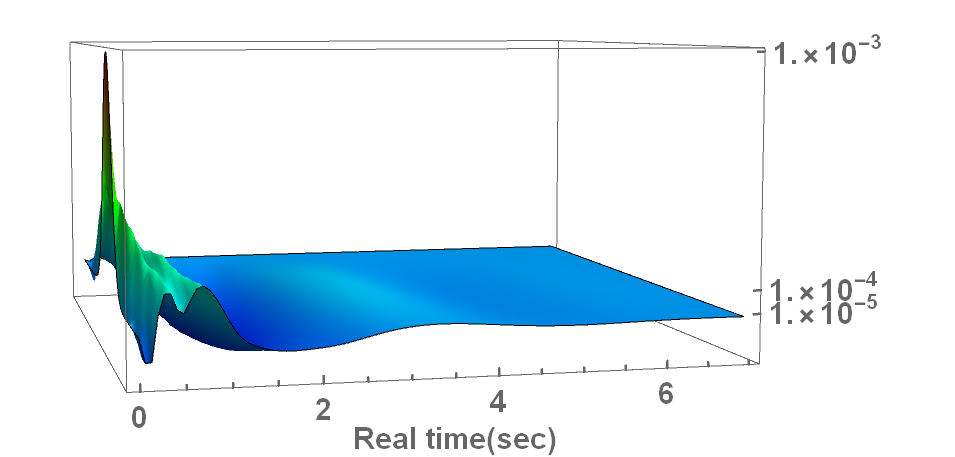} & \includegraphics[width=2.1in]{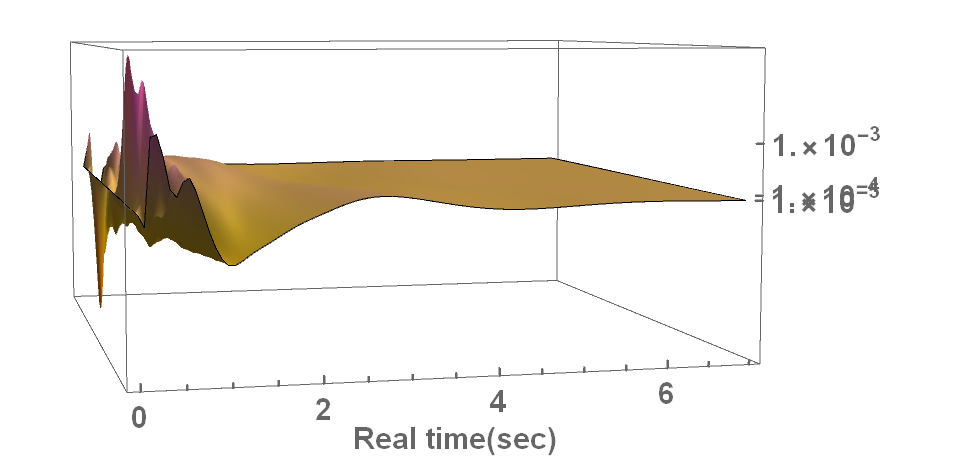} \\ \hline
\end{tabular}
\caption{\small Graphs in each column show simulations for the approximated solutions $\{v, w,\psi\}$ of the (M-T) closed-loop system for the uncontrolled, controlled, and partially controlled ($V(t)=0$) cases, respectively.  The $t-$axis is real-time.}
\label{g2}
\end{table}
%Here is the plot of the voltage controller feedback in time:
%\begin{table}[ht]
%\begin{tabular}{ll} \hline
%\centering Euler-Bernoulli (E-B) & Mindlin-Timoshenko (M-T) \\ \hline\\
%\includegraphics[width=3.0in]{EB-nl-un-v.png} & \includegraphics[width=3.0in]{M-T-controllers-vol.png} \\
%\includegraphics[width=3.0in]{EB-nl-dam-v.png} & \includegraphics[width=3.0in]{M-T-controllers-m.png} \\
%\includegraphics[width=3.0in]{EB-nl-dam-v-justbendingcont.png} & \includegraphics[width=3.0in]{M-T-controllers-g.png}  \\ \hline
%\end{tabular}
%\caption{\small Controllers in time.}
%\label{g3}
%\end{table}
 \bgroup
\def\arraystretch{1.1}
\begin{table}[htb]
\centering
\begin{tabular}{cc}
%(E-B) Total energy   & (M-T) Total energy \\
\includegraphics[width=3.0in]{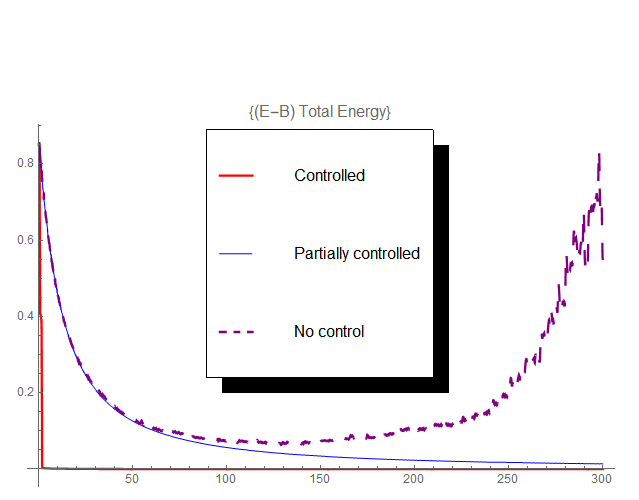} & \includegraphics[width=3.0in]{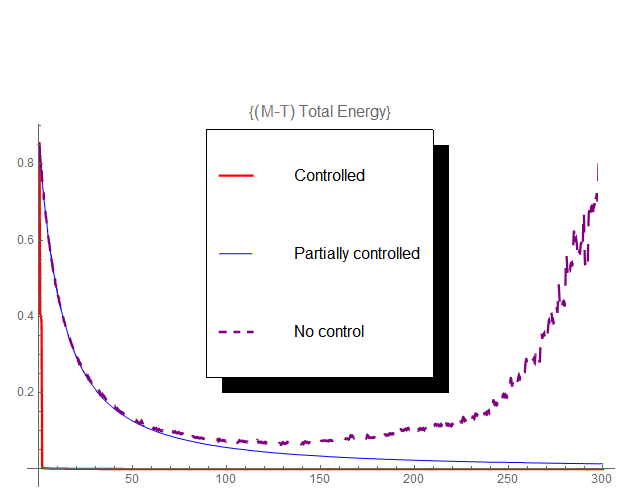}  
\end{tabular}
\caption{\footnotesize Total approximated energy  of the  (E-B)  and (M-T) controlled, partially controlled, uncontrolled  models, respectively. The $t-$axis is non-dimensional with the non-dimensionality constant $A_1=0.023$. As $dx\to 0$ (larger $N$), the actual decay can be observed better.}
\label{g3}
\end{table}

 \bgroup
\def\arraystretch{0.9}
\begin{table}[htb]
\begin{tabular}{ll}
%\centering (E-B) Tip velocities  & (M-T) Tip velocities \\
\includegraphics[width=2.8in]{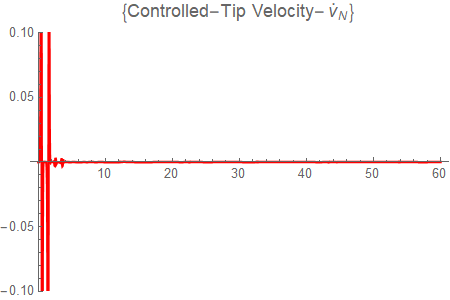} & \includegraphics[width=2.8in]{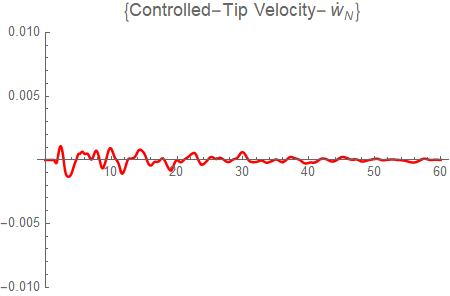}  \\
\includegraphics[width=2.8in]{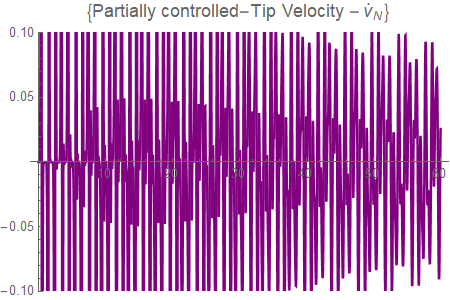} & \includegraphics[width=2.8in]{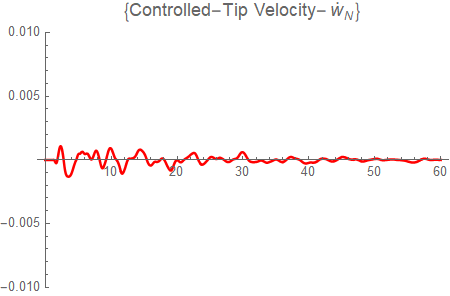}  \\
\includegraphics[width=2.8in]{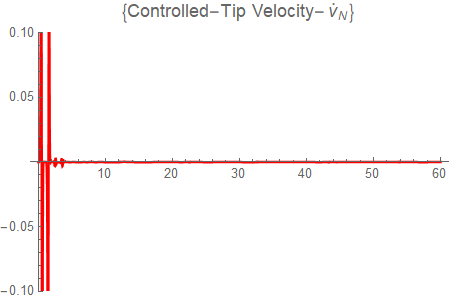} & \includegraphics[width=2.8in]{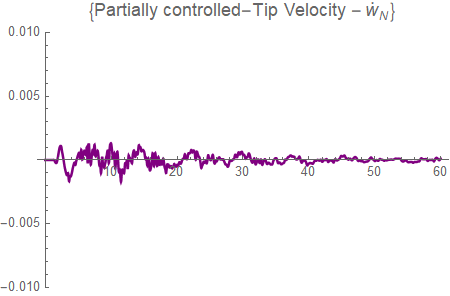}
\end{tabular}
\caption{\footnotesize Graphs in each column show simulations for the tip velocities $\{\dot v_N, \dot w_N\}$ of the (E-B)  closed-loop system for the controlled, partially controlled, and uncontrolled models. The $t-$axis is non-dimensional with the non-dimensionality constant $A_1=0.023$.  }
\label{g4}
\end{table}
 \bgroup
\def\arraystretch{0.9}
\begin{table}[htb]
\begin{tabular}{lll}
%\centering (E-B) Tip velocities  & (M-T) Tip velocities \\
\includegraphics[width=2.0in]{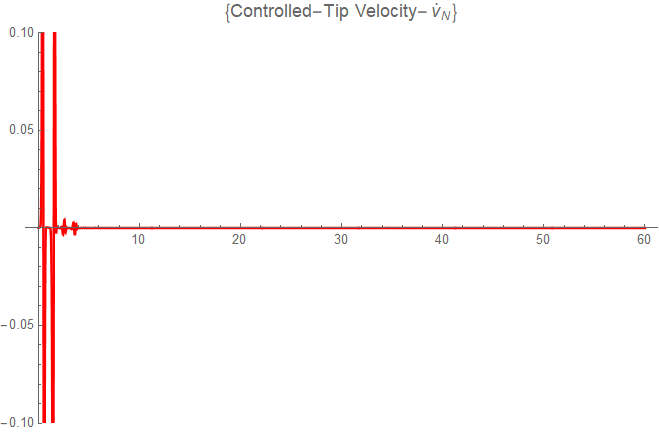} & \includegraphics[width=2.0in]{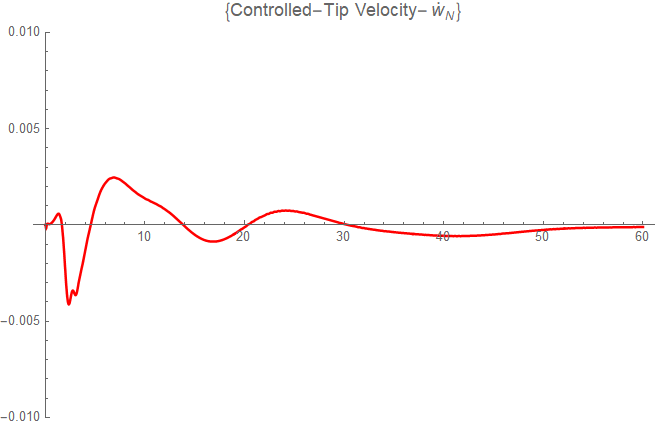}  &  \includegraphics[width=2.0in]{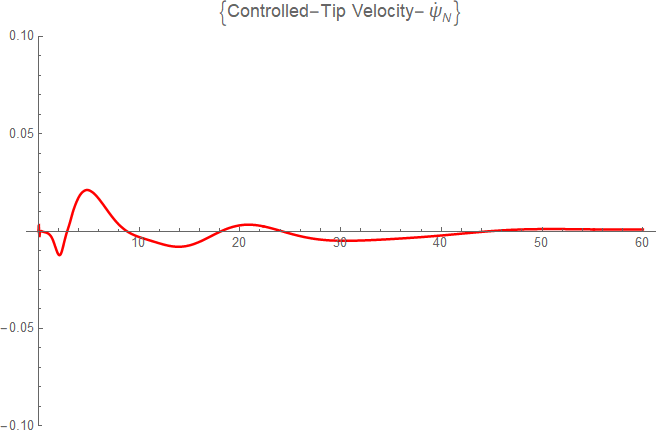} \\
\includegraphics[width=2.0in]{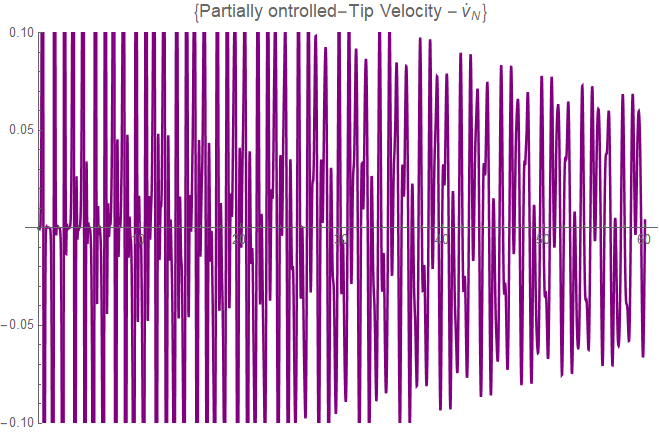} & \includegraphics[width=2.0in]{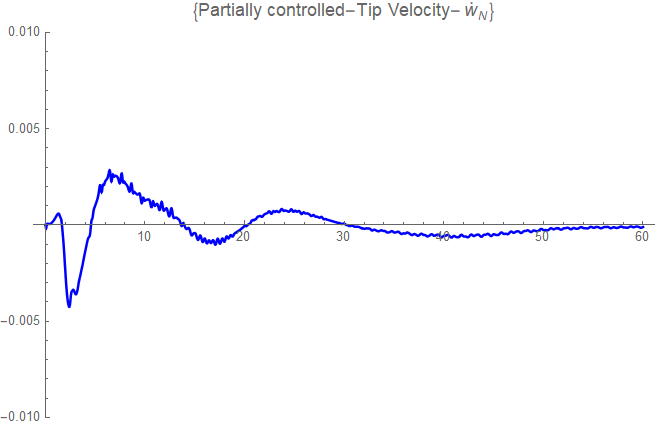} &  \includegraphics[width=2.0in]{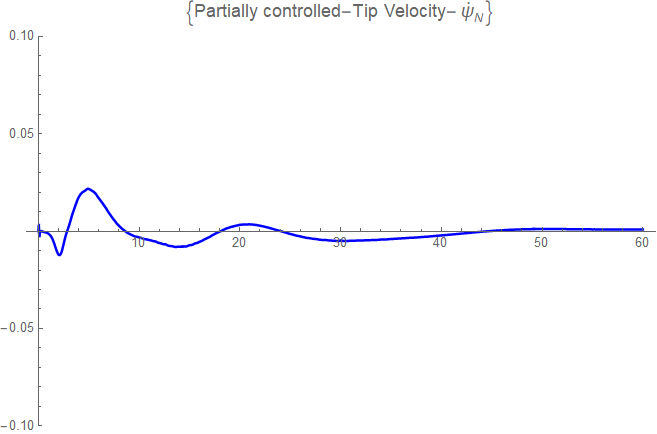} \\
\includegraphics[width=2.0in]{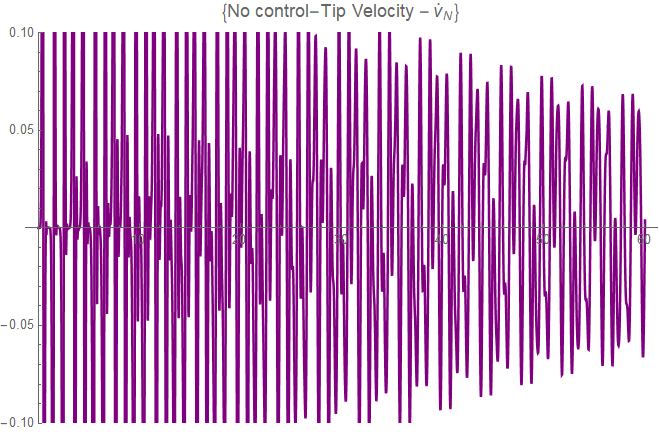}      & \includegraphics[width=2.0in]{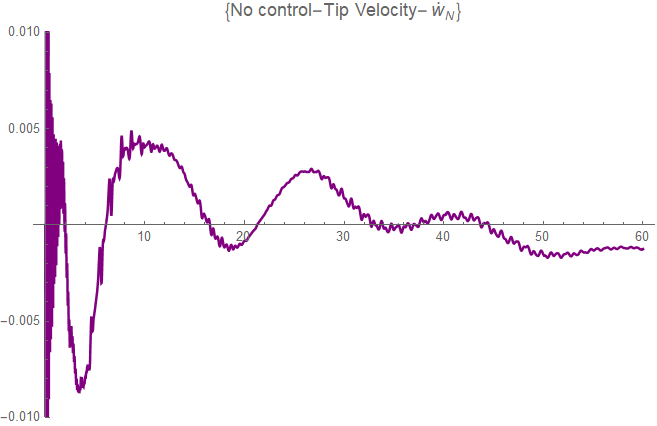} &  \includegraphics[width=2.0in]{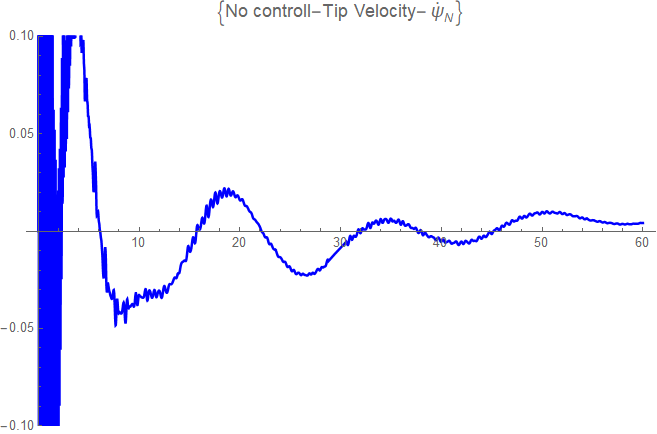} \\
\end{tabular}
\caption{\footnotesize Graphs in each column show simulations for  the tip velocities $\{\dot v_N, \dot \psi_N, \dot w_N\}$ of the (M-T) closed-loop system for the controlled, partially controlled, and uncontrolled models. The $t-$axis is non-dimensional with the non-dimensionality constant $A_1=0.023$.}
\label{g5}
\end{table}

\section{Discussion and Conclusion}
 A  large class of nonlinear piezoelectric smart beam models are obtained in the framework of a consistent variational approach. All models derived in this manuscript fit in the category of quasi-linear hyperbolic systems in \cite{Li}, and therefore, the well-posedness of these models can be shown in appropriate Hilbert spaces. Since the spectral analysis for these models is very tedious, the Finite Difference Method is very suitable for stable approximations.  Comparing our results obtained by the filtered Finite Differences with the ones  obtained by the mixed-Finite Element  method   is essential  \cite{Ito}. This together with the  extended numerical analysis with fourth order approximations  is currently under investigation  \cite{K-O}. It is also an ongoing project to find out the decay rate of the global stabilization in discrete/continuous of electrostatic and fully dynamic models with $B^*-$type feedback controllers \cite{Ozkan7}.

Another ongoing work is the implementation of the Lyapunov's direct method  to design non-trivial boundary controllers for the electrostatic and fully dynamic models. The designed controller's ability to stabilize the composite at its equilibrium position is proven analytically.%Fully dynamic and electrostatic models derived in this manuscript fit in the category of quasi-linear hyperbolic systems in \cite{Li}. The main challenge in our model  that the voltage control $V(t)$ appears at both  the axial strain boundary condition at $x=L$ and  the $w-$equation. The boundary controlling nature of $V(t)$ makes the bilinear control operator unbounded. Similar bilinear frameworks are recently studied in \cite{Ay,Ber}.
%$ yet it does not fit in the bilinear Euler-Bernoulli model studied in \cite{Lenhart} where the control operator is bounded.
%(i) The boundary controllability and stabilization of a nonlinear beam model (without considering the piezoelectricity) is studied in i.e. \cite{G-Li,Lag,Lag-L}, and the references therein. In fact,

For energy harvesting applications, three-layer models obtained in Section \ref{S3} are attached to a circuit equation, see Remark \ref{har}. It is a future work to carry out a detailed spectral analysis for the corresponding linearized models since it may be helpful for deducing stability characteristics of the nonlinear models rigorously. This analysis is similar to the one carried out for a unimorph  \cite{Shubov}.

%The well-posedness, controllability, and stabilizability problems are handled in details in \cite{O-M1}.

\subsection{Acknowledgments}
 This  research was supported by the Western Kentucky University startup grant. The author acknowledges the fruitful discussions with  Prof. Mikhail Khenner  on the Finite Difference approximations.

% References
%\bibliography{report} % bibliography data in report.bib

\begin{thebibliography}{99}
\bibitem{Wehbe} F. Abdallah, S. Nicaise, J. Valein, A. Wehbe, {Uniformly exponentially or polynomially
stable approximations for second order evolution equations and some applications,} {\sl ESAIM:
Control, Optimisation and Calculus of Variations} {\bf 19-3}, 844--887 (2013).
\bibitem{Alabou} F. Alabau-Boussouira, Y. Privat, E. Trélat. {Nonlinear damped partial differential
equations and their uniform discretizations,} {\sl Journal of Functional Analysis} {\bf 273-1}, 352--403 (2017).
\bibitem{A-C} K. Ammari, S. Nicaise, {Numerical approximations for the best decay rate for some dissipative systems,}  arXiv:1707.09155v1 [math.NA] (28 Jul 2017).
\bibitem{Ashgari} M. Asghari, M.H. Kahrobaiyan, M.T. Ahmadian, { A nonlinear Timoshenko beam formulation based
on the modified couple stress theory,} {\sl International Journal of Engineering Science} {\bf 48}, 1749--1761 (2010).
\bibitem{Ay} R. El Ayadi , M. Ouzahra,  A. Boutoulout, { Strong stabilisation and decay estimate for
unbounded bilinear systems, } {\sl International Journal of Control} {\bf 85-10} 1497-1505, (2012).
%\bibitem{Banks} H.T. Banks, R.C. Smith, Y. Wang, [Smart material structures: Modelling, Estimation and Control], Mason, Paris (1996).
\bibitem{Ito}H.T. Banks, K. Ito and B. Wang, Exponentially stable approximations of weakly damped wave equations,
{\sl International Series of Numerical Mathematics} {\bf 100}, 1--33 (1991).


\bibitem{Baz} A. Baz, {Boundary Control of Beams Using Active Constrained Layer Damping,} {\sl J. Vib. Acoust.} {\bf 119-2}, 166-172 (1997).
\bibitem{Bea} K. Beauchard, Local controllability of a 1-D beam equation, {\sl SIAM J. Control Optim.} {\bf 47-3}, 1219--1273 (2008).
\bibitem{Ber} L. Berrahmoune, {Stabilization of unbounded bilinear control systems in Hilbert space}, {\sl Journal of Mathematical Analysis and
Applications} {\bf 372}, 645--655  (2010).
%\bibitem{C-Z} R. F. Curtain and H. Zwart, {\sl  An Introduction to Infinite-dimensional
%Linear Systems Theory},  Volume 21 of Texts in Applied Mathematics, Springer-Verlag, New York, 1995.

\bibitem{Roventa} I. F. Bugariu, S. Micu, and I. Roventa, {Approximation of the controls for the beam equation with vanishing viscosity,} {\sl Mathematics of Computation} {\bf 85-11}, 2259--2303 (2016).
\bibitem{Cao-Chen} Y. Cao, X.B. Chen, {A Survey of Modeling and Control Issues for Piezo-electric Actuators,} {\sl Journal of Dynamic Systems, Measurement, and Control } {\bf 137-1}, 014001 (2014).
\bibitem{Dag} C. Dagdeviren, et al., {\sl Conformal piezoelectric energy harvesting and storage from motions of the heart, lung, and diaphragm,} Proc. Natl. Acad. Sci. U. S. A. {\bf 111}, 1927-1932 (2014).
\bibitem{Dag1} C. Dagdeviren, et al., {\sl Recent Progress in Flexible and Stretchable Piezoelectric Devices for Mechanical Energy Harvesting, Sensing and Actuation}, Extreme Mechanics Letter {\bf 9(1)}, 269-281  (2016).
\bibitem{Review} S Devasia, E. Eleftheriou, S. O. Reza Moheimani, {A Survey of Control Issues in Nanopositioning}, {\sl  IEEE Transactions on Control Systems Technology} {\bf 15-5},  802--823  (2007).
\bibitem{Dietl} J.M. Dietl, A.M. Wickenheiser, E. Garcia, {A Timoshenko beam model for cantilevered piezoelectric energy harvesters,} {\sl{Smart Mater. Struct.}} {\bf  19}, 055018 (2010).
\bibitem{E-Inman} A. Erturk, D. Inman, {A distrbiuted parameter model for cantilever model for piezoelectric energy harvesting from base excitations,} {\sl{J. Vib. Acoust.}} {\bf  130}, 041002 (2008).
 \bibitem{Exacto} {\sl EXACTO Guided Bullet Demonstrates Repeatable Performance against Moving Targets}, DARPA, 27 April 2015. Last accessed: 5 January 2018.
 \bibitem{health}  V. Giurgiutiu, [Structural Health Monitoring/with Piezoelectric Wafer Active Sensors], Elsevier Academic Press (2008).
\bibitem{Fur} K. Furutani, M. Urushibata, N. Mohri, Displacement control of piezoelectric element by feed-
back of induced charge, {\sl Nanotechnology} {\bf 9 },  93-98 (1998).
\bibitem{G-Li} Q. Gu, G. Leughering, T. Li,  {Exact Boundary Controllability on a Tree-Like Network of Nonlinear Planar Timoshenko Beams,} {\sl Chin. Ann. Math,  Series B} {\bf 38B(3)}, 711--740 (2017).

\bibitem{Gu} G.Y. Gu,  L.M. Zhu, C.Y. Su, H. Ding, S, Fatikow,  {Modeling and Control of Piezo-Actuated Nanopositioning Stages: A Survey,} {\sl  IEEE Transactions on Automation Science and Engineering} {\bf 13-1},     313--332 (2016).
\bibitem{Hansen3}  S.W. Hansen,
{ Several Related Models for Multilayer Sandwich Plates,}
{ \sl Mathematical Models \& Methods in Applied
Sciences}  {\bf 14-8}, 1103-1132  (2004).
%\bibitem{O-Hansen1}{S.W. Hansen, A.\"{O}. \"{O}zer}, \newblock{ Exact boundary controllability of an abstract Mead-Marcus Sandwich beam model,} \newblock{\sl The  Proceedings of $49^{\rm th}$ IEEE Conf. on Decision \& Control}, Atlanta, USA, (2010) 2578-2583.
%\bibitem{K-M-M2} B. Kapitonov, B. Miara, and G. P. Menzala, Stabilization of a layered 3--D body by boundary dissipation, {\sl{ESAIM:COCV}}, vol. 12, 2006, pp. 198--215.


\bibitem{Kugi} A. Kugi, S. Kurt, H. Irschik, {{Infinite-dimensional control of nonlinear beam vibrations by piezoelectric actuator and sensor layers,}} {\sl Nonlinear Dynamics} {\bf 19-1},   71--91 (1999).
\bibitem{Leon} L. Leon, E. Zuazua, {Boundary controllability of the finite-difference space semi-discretizations of the
beam equation,} {\sl ESAIM Control Optim. Calc. Var.} {\bf 8}, 827-862 (2002).
%\bibitem{Lag} J.E. Lagnese, {\sl{Boundary controllability of nonlinear beams to bounded states}}, Proc. Royal Soc. of Edinburg, {\bf 119A} (1991),  63--72.
\bibitem{Lagnese-Lions} J.E. Lagnese, J.-L. Lions, [Modeling Analysis and Control of Thin Plates], Masson, Paris (1988).
\bibitem{Lag-L} J.E. Lagnese, G. Leugering, {{Uniform stabilization of a nonlinear beam by nonlinear boundary feedback}}, {\sl J. Diff. Eqns} {\bf 91},  355-388 (1991).
\bibitem{Lam} M.J. Lam, D. Inman, W. R. Saunders, {Vibration Control through Passive Constrained Layer Damping and Active Control,}  {\sl J. Intel. Mater. Syst. Str.} {\bf 8-8},  663-677 (1997).


%\bibitem{L-M} I. Lasiecka and B. Miara, Exact controllability of a 3D piezoelectric body, {{\sl C. R. Math. Acad. Sci. Paris}}, {\bf  347} (2009),  167--172.
  %  \bibitem{Lenhart} { Bilinear spatial control of the velocity term in a
%Kirchhoff plate equation, }{\sl Electric Journal of Differential Equations,} {\bf 27} (2001) 1--15.
%\bibitem{Lee} P.C.Y. Lee,  A variational principle for the equations of piezoelectromagnetism in elastic dielectric crystals,{\sl{ Journal of Applied Physics,}} vol 69 (11), (1991), pp. 7470--7473.	
\bibitem{Li} Li, T. T., [Controllability and Observability for Quasilinear Hyperbolic Systems,] {\sl AIMS Ser. Appl. Math.} {\bf 3} American Institute of Mathematical Sciences and Higher Education Press (2010).
\bibitem{M-Z} A. Marica, E. Zuazua,  {[On the Lack of Uniform Observability for Discontinuous Galerkin Approximations of Waves.]} {In: Symmetric Discontinuous Galerkin Methods for 1-D Waves}, New York, NY, Springer (2014).
 \bibitem{Mead} D.J. Mead and S. Markus, {The forced vibration of a three-layer, damped sandwich beam
with arbitrary boundary conditions,} {\sl J. Sound Vibr.} {\bf 10}, 163--175  (1969).

\bibitem{O-M} K.A. Morris and  A.\"{O}. \"{O}zer, Strong stabilization of piezoelectric beams with magnetic effects, {\sl{The Proceeedings of the $52^{\text{nd}}$ Conference on Decision and Control,}} Florence, Italy, 3014--3019 (2014).
\bibitem{O-M1}  K.A. Morris, A.\"{O}. \"{O}zer, { Modeling and stabilizability of voltage-actuated piezoelectric beams with magnetic effects,}  {\sl SIAM J. Control Optim.} {\bf 52--4},  2371--2398 (2014).
 \bibitem{N-Z}   M. Negreanu, E. Zuazua, {Convergence of a multigrid method for the controllability of a 1-d wave equation,} {\sl C. R. Math. Acad. Sci. Paris} {\bf 338}, 413-418 (2004).
\bibitem{Ozkan2} A.\"{O}. \"{O}zer, {\newblock{Semigroup well-posedness of a voltage controlled active constrained layered (ACL) beam with magnetic effects,}} {\sl The Proceedings of the American Control Conference,} Boston, MA, USA, 4580-4585 (2016).
    \bibitem{Ozkan3} A.\"{O}. \"{O}zer, {\newblock Modeling and controlling an active constrained layered (ACL) beam actuated by two voltage sources with/without magnetic effects, } {\sl IEEE Transactions of Automatic Control}, {\bf 62-12},  6445--6450  (2017).
    \bibitem{Ozkan6}  {A.\"O. \"Ozer,}   {\sl Exponential stabilization of the smart piezoelectric composite beam with only one
boundary controller, } {\sl The Proceedings of the  International Federation of Automatic Control Conference}, in press (2018).
        \bibitem{Ozkan7}  {A.\"O. \"Ozer,}   {Well-posedness and stabilization analysis of nonlinear piezoelectric devices, } in preparation.
\bibitem{O-Hansen3} {A.\"{O}. \"{O}zer, and S.W. Hansen,} {Exact boundary controllability results for a multilayer Rao-Nakra sandwich beam}, {\sl SIAM J. Control Optim.,} {\bf 52-2}, 1314--1337 (2014).

\bibitem{O-Hansen4} {A.\"O. \"Ozer, S.W. Hansen}, {Uniform stabilization of a multi-layer Rao-Nakra sandwich beam},  {\sl Evolution Equations and Control Theory} {\bf 2-4}, 195--210 (2013).
\bibitem{K-O}   A.\"O. \"Ozer, M. Khenner, {Numerical investigation of  boundary controlled nonlinear piezoelectric beams,} in preparation.
\bibitem {Rao} Y.V.K.S. Rao and B.C. Nakra, {Vibrations of unsymmetrical sandwich beams and plates with viscoelastic cores,} {\sl J. Sound Vibr.} {\bf 34-3}, 309--326  (1974).
  \bibitem{Ru} C. Ru, X. Liu, Y. Sun, [Nanopositioning Technologies: Fundamentals and Applications], Springer International  (2016).
%\bibitem{Scott} W. T. Scott, {\sl{Approximation to real irrationals by certain classes of rational fractions,}} Bull. Amer. Math. Soc. vol. 46 (1940) pp. 124-129.
\bibitem{B-R} M.C. Ray, A. Baz, {Control of Nonlinear Vibration of Beams Using Active Constrained Layer Damping}, {\sl Journal of Vibration and Control} {\bf 7-4}, 539--549 (2001).

\bibitem{Ronkanen} P. Ronkanen, P. Kallio, M. Vilkko, H.N. Koivo,  Displacement Control of Piezoelectric Actuators Using Current and Voltage, {\sl IEEE/ASME Trans. Mechatronics} {\bf 16-1}, 160-166 (2011).
%\bibitem{Rogacheva} N. Rogacheva, {\sl{ The Theory of Piezoelectric Shells and Plates}}, Boca Raton, FL: CRC Press; 1994.
%\bibitem{Sene} A. Sen\`{e}, {\sl{Modelling of piezoelectric static thin plates   ,}} { Asymptotic Analysis},
%(25-1) (2001), pp. 1--20.
 \bibitem{Sci-D} {\sl Ultrasound imaging of the brain and liver,} Science Daily Magazine, Source: Acoustical Society of America, 26 June 2017. Last Accessed: 5 January 2018.
\bibitem{Shen}  I.Y. Shen, {A variational formulation, a work-energy relation and damping mechanisms of active constrained
layer treatments,} {\sl Journal of Vibration and Acoustics} {\bf 119-2}, 192-199, (1997).
\bibitem{Shubov} M.A. Shubov, {Spectral analysis of a non-selfadjoint operator generated by an energy harvesting model and application to an exact controllability problem}, Asymptotic Analysis {\bf 102(3-4)}, 119-156 (2017).
\bibitem{Smith} R.C. Smith,  [Smart Material Systems], Society for
Industrial and Applied Mathematics (2005).
    \bibitem{T-Z} L.T. Tebou, E. Zuazua. {Uniform boundary stabilization of the finite difference space discretization of
the 1-d wave equation.} {\sl Adv. Comput. Math.} {\bf 26}, 337-365 (2007).
\bibitem{Tiersten} H.F. Tiersten, [Linear Piezoelectric Plate Vibrations], Plenum Press, New York (1969).
%\bibitem{Trigg0} \newblock R. Triggiani, {\newblock\emph{On the stabilizability  problem in Banach space,}} \newblock {J. Math. Anal. Appl.} (3), \textbf{52}  (1975), 383--403.
%\bibitem{Tucs} M. Tucsnak, Regularity and exact controllability for a beam with piezoelectric
%actuators, {\sl{ SIAM J. Cont. Optim.,}} vol. 34, 1996, pp. 922--930.
%\bibitem{Weiss-Tucsnak} M. Tucsnak and G. Weiss, {\sl{Observation and Control for Operator Semigroups}}, Birkhäuser Verlag, Basel (2009).
\bibitem {Ultra} {\sl Ultrasonic Welding: Basic Knowledge,} Beijing Ultrasonic [PDF document], uploaded on 07/01/2013.
%\bibitem{Weiss} M. Tucsnak, G. Weiss,  {\sl{How to get a conservative well-posed linear system out of thin air, Part II: controllability and stability,}} SIAM J. Cont. Optim. (42-3) (2003), pp. 907-935.
\bibitem{Voss} T. Voss, J. M. A. Scherpen, {Port-Hamiltonian Modeling of a Nonlinear Timoshenko Beam with Piezo Actuation,} {\sl SIAM J. Control Optim.,} {\bf 52-1}, 493--519 (2014).
%\bibitem{Tzou} H.S. Tzou, {\sl{Piezoelectric shells, Solid Mechanics and Its applications 19,}} Kluwer Academic, The Netherlands; 1993.
%\bibitem{Wang-Guo} J-M Wang, B-Z Guo {\sl On the stability of swelling porous elastic soils with fluid
%saturation by one internal damping,} IMA Journal of Applied Mathematics (71)(2006), pp. 565-582.
%\bibitem{Young} R.M. Young, {\sl{An Introduction to Nonharmonic Fourier Series,}} New York, Academic Press, (1980).
    \bibitem{Yang} J. Yang, [Special Topics in the Theory of Piezoelectricity], Springer, New York  (2009).
    %\bibitem{Yang1} J. Yang, A review of a few topics in piezoelectricity, {\sl{ Appl. Mech. Rev.,}} vol. 59, 2006, pp. 335–-345.
%\bibitem{Zhang} C-G Zhang, Regularity and exact controllability for the Timoshenko beam with Piezoelectric actuator,  Rocky Mountain J. Math., vol. 41 (3), 2011, pp. 999--1010.
\end{thebibliography}
%\bibliographystyle{spiebib} % makes bibtex use spiebib.bst

\end{document}